\documentclass[10pt,reqno]{extarticle}

\usepackage{geometry}
\usepackage{graphicx}

\usepackage{tocloft}
\usepackage{titlesec}
\titleformat{\section}{\normalfont\rmfamily\Large\bfseries\color{black}}{\thesection}{0.8em}{}
\titleformat{\subsection}{\normalfont\rmfamily\large\bfseries\color{black}}{\thesubsection}{0.8em}{}
\titleformat{\subsubsection}{\normalfont\rmfamily\normalsize\bfseries\color{black}}{\thesubsubsection}{0.8em}{}

\usepackage[utf8]{inputenc}
\usepackage[english]{babel}

\usepackage{amssymb}
\usepackage{amsmath}
\usepackage{amsthm}
\usepackage{mathrsfs}
\usepackage{mathtools}

\usepackage[all,arc]{xy}
\usepackage{tikz}
\usepackage{tikz-cd}
\usepackage{float}
\usepackage{caption}
\usepackage{subcaption}

\usepackage{soul}
\usepackage{cancel}
\usepackage{slashed}
\usepackage{tensor}
\usepackage{scalerel}

\usepackage{color}
\usepackage{enumitem}
\usepackage{hyperref}
\hypersetup{
  colorlinks   = true, 
  urlcolor     = blue, 
  linkcolor    = blue, 
  citecolor   = red 
}

\usepackage[final]{pdfpages}
\usepackage{showlabels}

\numberwithin{equation}{section}
\allowdisplaybreaks

\newtheorem{theorem}{Theorem}[section]
\newtheorem{definition}[theorem]{Definition}
\newtheorem{corollary}[theorem]{Corollary}
\newtheorem{remark}[theorem]{Remark}
\newtheorem{lemma}[theorem]{Lemma}
\newtheorem{proposition}[theorem]{Proposition}
\newtheorem{problem}[theorem]{Problem}
\newtheorem{claim}[theorem]{Claim}

\newcommand{\ed}{\mathrm{d} \makebox[0ex]{}}
\newcommand{\calF}{\mathcal{F} \makebox[0ex]{}}
\newcommand{\circg}{\overset{\raisebox{-0.6ex}[0.6ex][0ex]{$\scaleto{\circ}{0.6ex}$}}{g} \makebox[0ex]{}}
\newcommand{\uC}{\underline{C} \makebox[0ex]{}}
\newcommand{\uu}{\underline{u} \makebox[0ex]{}}
\newcommand{\slashg}{g \mkern-8.2mu \scaleto{\boldsymbol{\slash}}{1.6ex} \mkern+1mu \makebox[0ex]{}}
\newcommand{\subslashg}{g \mkern-7.7mu \scaleto{\boldsymbol{\slash}}{1ex} \mkern+1mu \makebox[0ex]{}}

\newcommand{\uL}{\underline{L} \makebox[0ex]{}}
\newcommand{\uchi}{\underline{\chi} \makebox[0ex]{}}
\newcommand{\tr}{\mathrm{tr} \makebox[0ex]{}}
\newcommand{\hatchi}{\hat{\chi} \makebox[0ex]{}}
\newcommand{\hatuchi}{\underline{\hat{\chi}} \makebox[0ex]{}}
\newcommand{\circnabla}{\overset{\raisebox{0ex}[0.6ex][0ex]{$\scaleto{\circ}{0.6ex}$}}{\raisebox{0ex}[1.2ex][0ex]{$\nabla$}} \makebox[0ex]{}}
\newcommand{\dpartial}{\dot{\partial} \makebox[0ex]{}}

\newcommand{\circDelta}{\overset{\raisebox{0ex}[0.6ex][0ex]{$\scaleto{\circ}{0.6ex}$}}{\raisebox{0ex}[1.2ex][0ex]{$\Delta$}} \makebox[0ex]{}}
\newcommand{\R}{\mathrm{R} \makebox[0ex]{}}

\newcommand{\dvol}{\mathrm{dvol} \makebox[0ex]{}}
\newcommand{\calQ}{\mathcal{Q} \makebox[0ex]{}}

\newcommand{\bt}{\bar{t} \makebox[0ex]{}}
\newcommand{\br}{\bar{r} \makebox[0ex]{}}
\newcommand{\btheta}{\bar{\theta} \makebox[0ex]{}}
\newcommand{\buu}{\underline{\bar{u}} \makebox[0ex]{}}
\newcommand{\barf}{\bar{f} \makebox[0ex]{}}
\newcommand{\bh}{\bar{h} \makebox[0ex]{}}
\newcommand{\bv}{\bar{v} \makebox[0ex]{}}

\newcommand{\calL}{\mathcal{L} \makebox[0ex]{}}

\newcommand{\balpha}{\bar{\alpha} \makebox[0ex]{}}

\newcommand{\bvartheta}{\bar{\vartheta} \makebox[0ex]{}}

\title{\textsc{A Lorentz invariant sharp Sobolev inequality on the circle}}
\author{Pengyu Le}
\newcommand{\Address}{{
  \bigskip
  \footnotesize
  \textsc{Yanqi Lake Beijing Institute of Mathematical Sciences and Applications, Beijing, China}
  
  \textit{E-mail address}: \texttt{pengyu.le@bimsa.cn}
}}
\date{}

\begin{document}

\maketitle

\begin{abstract}
We prove the following sharp Sobolev inequality on the circle
$$\textstyle \int_{\mathbb{S}^1} [4(v')^2 - v^2]  \ed \theta \geq  - \frac{4\pi^2}{\int_{\mathbb{S}^1} v^{-2} \ed \theta},$$
with the equality being achieved when
$v^{-2} (\theta) = \frac{k\sqrt{1-\alpha^2}}{1+ \alpha \cos(\theta - \theta_o)}$
where
$k>0$, $\alpha \in (-1,1)$, $\theta_0 \in \mathbb{R}$.
If $v$ vanishes somewhere on the circle, then
$$\textstyle 4 \int_{\mathbb{S}^1} (v')^2 \ed \theta \geq\int_{\mathbb{S}^1} v^2  \ed \theta.$$
The basic tools to prove the inequality are the rearrangement inequality on the circle and the variational method. We investigate the variational problem of the functional $\calF[v] = \int_{\mathbb{S}^1} [4(v')^2 - v^2] \ed \theta$ under the constraint $\int_{\mathbb{S}^1} v^{-2} \ed \theta = 2\pi$. 
An important geometric insight of the functional $\calF$ is that it is invariant under the Lorentz group, since $\calF[v]$ is the integral of the product of two null expansions of a spacelike curve parameterised by the function $v^{-2}$ in a lightcone in $3$-dim Minkowski spacetime. The global minimiser of $\calF$ under the constraint is simply given by the spacelike plane section of the lightcone.
We introduce a method which combines the symmetric decreasing rearrangement and the Lorentz transformation. This method isnot confined to the scope of this paper, but is applicable to other Lorentz invariant variational problems on $\mathbb{S}^{n}, n \geq 1$. As an example, we sketch a proof of the sharp Sobolev inequality on $\mathbb{S}^n, n\geq 3$ by this method.

\end{abstract}

\tableofcontents

\section{Introduction}\label{sec 1}

In this paper, we prove a sharp Sobolev inequality on the circle that
\begin{align}
\int_{\mathbb{S}^1} [4(v')^2 - v^2]  \ed \theta \geq  - \frac{4\pi^2}{\int_{\mathbb{S}^1} v^{-2} \ed \theta},
\label{eqn 1.1}
\end{align}
with the equality being achieved when
$v^{-2} (\theta) = \frac{k\sqrt{1-\alpha^2}}{1+ \alpha \cos(\theta - \theta_o)}$
where
$k>0$, $\alpha \in (-1,1)$, $\theta_0 \in \mathbb{R}$. See theorem \ref{thm 7.1}. As a corollary of the above, we show that if $v$ vanishes somewhere on the circle, then
\begin{align*}
4 \int_{\mathbb{S}^1} (v')^2 \ed \theta \geq\int_{\mathbb{S}^1} v^2  \ed \theta.
\end{align*}

We approach the about inequality through a geometric consideration. We depart from investigating the geometry of a spacelike section $S$ of a lightcone in the Minkowski spacetime and in particular studying a functional $\calF[S]$ which is the integral of the product of two null expansions relative to a conjugated null frame of the spacelike section. This integral is closely related to the Hawking mass in general relativity, see \cite{H1968}. We observe the following interesting facts about this functional (see subsection \ref{subsec 2.3} for the details):
\begin{enumerate}[label=\textit{\alph*.}]
\item
when the dimension of the Minkowski spacetime is $\geq 5$, we observe an interesting connection between this functional and the Yamabe invariant of the conformal class of the standard round sphere;
\item 
when the dimension of the Minkowski spacetime is $4$, the functional is constant which is simply a multiple of $2\pi$, the Euler characteristic of the $2$-dimensional sphere;
\item
when the dimension of the Minkowski spacetime is $3$, the functional $\calF[S]$ takes the simple form
\begin{align*}
\calF[S] = \int_{\mathbb{S}^1} [4 (v')^2 - v^2] \ed \theta,
\end{align*}
where $v$ is related to the parametrisation of $S$ in the lightcone.
\end{enumerate}
In cases \textit{a.} and \textit{c.}, it is natural to consider the corresponding variational problem of $\calF[S]$ under the constraint that the area (length in case \textit{c.}) $|S|$ is a constant. The variational problem in case \textit{a.} is equivalent to the Yamabe problem for the standard round sphere, which was solved originally in \cite{A1976a} using essentially the optimal Sobolev inequality on $\mathbb{R}^n$ proved independently in \cite{A1976b} and \cite{Ta1976}. We shall mention that the Yamabe problem for the compact and not locally conformally flat Riemannian manifold with dimension $\geq 6$ is solved in \cite{A1976a}, while the remaining case of dimension $3,4,5$ or being locally conformally flat was solved in \cite{S1984} making use of the positive mass theorem proved in \cite{SY1979}\cite{SY1988}. In case \textit{c.}, the variational problem is stated as problem \ref{pro 2.4}, which can be also phrased as follows.
\begin{problem}\label{pro 1.1}
Find the extreme of the functional $\calF[v] = \int_{\mathbb{S}^1} [4 (v')^2 - v^2] \ed \theta$ under the constraint $\int_{\mathbb{S}^1} v^{-2} \ed \theta = 2\pi$. Moreover find the function achieving the extreme.
\end{problem}
Problem \ref{pro 1.1} is the main motivation of the paper and the Sobolev inequality \eqref{eqn 1.1} gives its answer. The Sobolev type inequalities on the sphere $\mathbb{S}^n$, $n\geq 2$ have been studied extensively by many works before. For example for the $2$-dimensional standard round sphere case, the Onofri inequality obtained in \cite{O1982} states that if $w \in \mathrm{H}^1(\mathbb{S}^2)$, then $v$ satisfies
\begin{align}
\frac{1}{4\pi} \int_{\mathbb{S}^2} |\nabla v|^2 \ed \mu 
+ \frac{ \int_{\mathbb{S}^2} 2 v \ed \mu}{4\pi}
\geq
\log \frac{\int_{\mathbb{S}^2} e^{2v} \ed \mu}{4\pi}.
\label{eqn 1.2}
\end{align}
The Onofri inequality is closely related to the Moser-Trudinger inequality in \cite{Tr1967}\cite{M1971}\cite{A1979}. For the higher dimensional round sphere $\mathbb{S}^n$, $n \geq 3$, the sharp Sobolev inequality says that
\begin{align}
\int_{\mathbb{S}^n} | \nabla v|^2 \ed \mu + \frac{n(n-2)}{4} \int_{\mathbb{S}^n} v^2 \ed \mu 
\geq
\frac{n(n-2)}{4} |\mathbb{S}^n|^{\frac{2}{n}} \Big( \int_{\mathbb{S}^n} v^{\frac{2n}{n-2}} \ed \mu \Big)^{\frac{n-2}{n}},
\label{eqn 1.3}
\end{align}
which is equivalent to the sharp Sobolev inequality on $\mathbb{R}^n$, $n\geq 3$ in \cite{A1976b}\cite{Ta1976} by the stereographic transformation. Note that formally let $n=1$ in the above inequality, it becomes the inequality \eqref{eqn 1.1}. It is interesting that when $n=1$, the right hand side involves the negative power of $v$.

The Onofri inequality \eqref{eqn 1.2} and the sharp Sobolev inequality \eqref{eqn 1.3} are both invariant under the conformal transformation of $\mathbb{S}^n$, $n\geq 2$, see \cite{O1982}\cite{L1983}\cite{CL1990}\cite{LL2001}\cite{C2004}. Inequality \eqref{eqn 1.1} has a similar feature as the above two inequalities, that inequality \eqref{eqn 1.1} is invariant under the Lorentz transformation because of the geometric meaning of the functional $\calF$ in $3$-dimensional Minkowski spacetime, see section \ref{sec 3} and proposition \ref{pro 3.1}. The invariance can be explicitly stated as follows. Let $\nu_{\alpha,\btheta_0}(\btheta) = \frac{\sqrt{1-\alpha^2}}{1- \alpha \cos (\btheta - \btheta_0)}$, and $\varphi_{\alpha, \theta_0, \btheta_0}$ be the transformation between $\mathbb{S}^1$ that
\begin{align*}
\varphi_{\alpha, \theta_0, \btheta_0}:
\quad
\mathbb{S}^1 \rightarrow \mathbb{S}^1,
\quad
\theta \mapsto \btheta,
\end{align*}
where
\begin{align*}
\cos  (\btheta - \btheta_0)
=
\frac{\alpha + \cos (\theta-\theta_0)}{1+ \alpha  \cos(\theta-\theta_0)},
\quad
\sin  (\btheta - \btheta_0)
=
\frac{  \sqrt{1-\alpha^2} \sin(\theta-\theta_0)}{ 1+ \alpha  \cos(\theta-\theta_0) },
\end{align*}
then define the transformation $\gamma_{\alpha, \theta_0, \btheta_0}$ for functions on the circle by that $\bv = \gamma_{\alpha, \theta_0, \btheta_0}(v)$ with 
\begin{align}
\bv(\btheta) = [\nu_{\alpha, \btheta_0}(\btheta)]^{-\frac{1}{2}} v(\theta),
\quad
\btheta = \varphi_{\alpha, \theta_0, \btheta_0}(\theta).
\label{eqn 1.4}
\end{align}
We have that
\begin{align*}
\calF[\bv] = \calF[v],
\quad
\int_{\mathbb{S}^1} \bv^{-2}(\btheta) \ed \btheta = \int_{\mathbb{S}^1} v^{-2}(\theta) \ed \theta.
\end{align*}
To obtain a clearer analogy between the Lorentz invariance of inequality \eqref{eqn 1.1} and the conformal invariance of the Onofri inequality \eqref{eqn 1.2} and the sharp Sobolev inequality \eqref{eqn 1.3}, we note that the conformal group of the standard round sphere $\mathbb{S}^n$ is isomorphic to the Lorentz group of the $(n+2)$-dimensional Minkowski space.

The role of the conformal invariance for proving inequalities \eqref{eqn 1.2} and \eqref{eqn 1.3} is two-sided: when applying the direct method of calculus of variations to prove the inequalities, the non-compactness of the conformal group of the round sphere contributes a difficulty, while the conformal invariance also gives a freedom to modify the minimising sequence of the corresponding functional such that the modified sequence could actually converges or weakly converges to a extreme point of the functional. This is exactly the case for the proofs of the Onofri inequality in \cite{O1982} and the proof of the sharp Sobolev inequality \eqref{eqn 1.3} in \cite{L1983}\cite{CL1990}. 

An important idea in the proofs of the Moser-Trudinger inequality, the Onofri inequality and the sharp Sobolev inequality \eqref{eqn 1.3} is the rearrangement inequality which says that the symmetric decreasing rearrangement preserves the $L^p$ norm and doesnot increase the $W^{1,p}$ norm. The rearrangement inequality reduces the above inequalities to $1$-dimensional problems on $\mathbb{R}_{\geq 0}$.

As mentioned above, the Onofri inequality and the sharp Sobolev inequality \eqref{eqn 1.3} have another nice feature of the conformal invariance. \cite{CL1990} exploited this conformal invariance to give a conceptual proof of the sharp Sobolev inequality \eqref{eqn 1.3}. It introduced the method of competing symmetries employing the symmetric decreasing rearrangement on $\mathbb{R}^n$ and the rotation of $\mathbb{S}^n$ (a conformal mapping on $\mathbb{R}^n$) alternatively to construct a sequence converging to the extreme function of the corresponding functional. The name ``competing symmetries" refers to the spherical symmetry by the symmetric decreasing rearrangement and the conformal symmetry of the limit function. See also \cite{LL2001}\cite{B2009} for the expositions of this method.

We prove the sharp Sobolev inequality \eqref{eqn 1.1} on the circle by a method using the symmetric decreasing rearrangement and the Lorentz transformation in another way. We give an overview of this method in the following. 

The basic idea is the following geometric observation for the Lorentz transformation $\gamma_{\alpha,0,0}$: if $v$ is positive continuous and doesnot attend its minimum on $[-\frac{\pi}{2}, \frac{\pi}{2}]$, then there exists $\alpha>0$ such that the transformation $\gamma_{\alpha,0,0}$ increases the minimum of $v$, i.e. $\min\{ \gamma_{\alpha,0,0} (v) \} > \min\{v\}$. See lemmas \ref{lem 6.1} and \ref{lem 6.2}. Then with this observation, we can apply the symmetric decreasing rearrangement and the Lorentz transformation $\gamma_{\alpha,0,0}$ alternatively to construct a sequence $\{v_n\}$ such that every $v_n$ is symmetric decreasing on $[-\pi,\pi]$ and
\begin{align*}
\int_{\mathbb{S}^1} v_n^{-2} \ed \theta = 2\pi,
\quad
\calF[v_n] \text{ is nonincreasing and } \min\{v_n\} \text{ is nondecreasing as } n \rightarrow +\infty.
\end{align*}
See definition \ref{def 6.4} for the construction. The sequence $\{v_n\}$ is bounded from below by definition, and moreover we can actually show that such constructed $\{v_n\}$ is also bounded from above. See lemma \ref{lem 6.11}. Investigating the convergence of the sequence $\{v_n\}$, the limit function of the sequence $\{v_n\}$ must be constant on $[\frac{\pi}{2}, \pi]$, see lemma \ref{lem 6.8}. Thus we introduce the following set of functions $S_c$ that $v\in S_c$ if
\begin{enumerate}[label=\alph*.]
\item $v$ is symmetric nonincreasing on $[-\pi,\pi]$,
\item $v \in \mathrm{H}^1(\mathbb{S}^1)$,
\item $v$ is constant on $[\frac{\pi}{2}, \pi]$,
\item $\int_0^{2\pi} v^{-2}\, \ed \theta = 2\pi$.
\end{enumerate}
We prove that the infimum of $\calF[v]$ under the constraint  $\int_0^{2\pi} v^{-2}\, \ed \theta = 2\pi$ is equal to the infimum of $\calF[v]$ in $S_c$ and it can be achieved if and only if the infimum can be achieved in $S_c$. See proposition \ref{pro 6.13}. The variational problem of $\calF[v]$ in the set $S_c$ is much simpler since
\begin{enumerate}
\item there is automatically a lower bound for a function $v\in S_c$,
\item it is essentially a variational problem on the half circle, which eliminates the invariance by the non-compact Lorentz group.
\end{enumerate}
Thus we successfully prove the inequality \eqref{eqn 1.1} by studying the variational problem of $\calF[v]$ in $S_c$ in section \ref{sec 7}.

We shall remark that the above method employed to prove inequality \eqref{eqn 1.1} could be useful for other Lorentz invariant variational problems. As an example, we sketch a proof of the sharp Sobolev inequality \eqref{eqn 1.3} on $\mathbb{S}^n$, $n\geq 3$ by this method in section \ref{sec 9}.

\section{Motivation of the problem}\label{sec 2}

In this section, we first give a presentation on the geometry of spacelike surfaces (curves if in the $3$-dim Minkowski spacetime) in a lightcone. Then we introduce the functional $\calF$ considered in this paper.

\subsection{Parameterisation of a spacelike surface in a lightcone}\label{subsec 2.1}
\begin{figure}[H]
\centering
\begin{subfigure}[b]{0.3\textwidth}
\centering
\begin{tikzpicture}[scale=0.5]
\draw[->] (-4,0)--(4,0) node[below] {$x^1$};
\draw[->] (0,-2)--(0,4) node[right] {$x^0$};
\draw[->] (-3,-1.5)--(3,1.5) node[right] {$x^n$};
\end{tikzpicture}
\caption{Rectangular coordinates $\{x^0, x^1, \cdots, x^n\}$.}
\end{subfigure}
\hfill
\begin{subfigure}[b]{0.3\textwidth}
\centering
\begin{tikzpicture}[scale=0.5]
\draw[->] (-4,0)--(4,0);
\draw[->] (0,-2)--(0,4) node[right] {$t$};
\draw[->] (-3,-1.5)--(3,1.5);
\draw (2,-0.1)
to [out=-150,in=-5] (-1,-0.5)
to [out=175,in=-150] (-2,0.1)
to [out=30,in=175] (1,0.5)
to [out=-5,in=30] (2,-0.1);
\draw (2*1.5,-0.1*1.5)
to [out=-150,in=-5] (-1*1.5,-0.5*1.5)
to [out=175,in=-150] (-2*1.5,0.1*1.5)
to [out=30,in=175] (1*1.5,0.5*1.5)
to [out=-5,in=30] (2*1.5,-0.1*1.5);
\draw[->] (1.5,-0.75) -- (3,-1.5) node[right] {$\partial_r$};
\end{tikzpicture}
\caption{Spatial polar coordinates $\{t,r, \vartheta \in \mathbb{S}^{n-1}\}$.}
\end{subfigure}
\hfill
\begin{subfigure}[b]{0.3\textwidth}
\centering
\begin{tikzpicture}[scale=0.5]
\draw[->] (0,0)--(0,4);
\node[below] at (0,0) {\tiny$(u,\uu) =(0,0)$}; 
\draw (-3,3) -- (0,0) -- (3,3) node[above] {\tiny $C_u$};
\draw (-3,-1) -- (0,2) -- (3,-1) node[below] {\tiny $\uC_{\uu}$};
\draw (0.8,0.8)
to [out=-135,in=-45] (-0.8,0.8)
to [out=135,in=-135] (-0.8,1.2)
to [out=45,in=135] (0.8,1.2)
to [out=-45,in=45] (0.8,0.8); 
\draw[->]  (3,0)--(4,1) node[above] {\tiny $\partial_{\uu}$};
\draw[->]  (3,0)--(2,1) node[above] {\tiny $\partial_u$};
\end{tikzpicture}
\caption{Double null coordinates $\{u, \uu, \vartheta \in \mathbb{S}^{n-1}\}$.}
\end{subfigure}
\caption{Coordinate systems of $\mathbb{M}^{n+1}$}
\label{fig 1}
\end{figure}
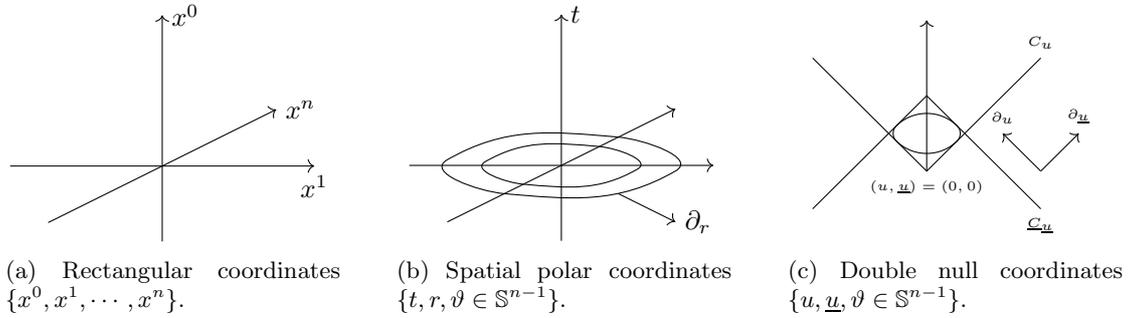

Let $(\mathbb{M}^{n+1},g)$ be the $(n+1)$-dim Minkowski spacetime. Consider three coordinate systems (see figure \ref{fig 1}):
\begin{enumerate}[label=\alph*.]
\item retangular coordinates: $\{x^0, x^1, \cdots, x^n\}\quad g = - (\ed x^0)^2 + (\ed x^1)^2 + \cdots + (\ed x^n)^2$;

\item spatial polar coordinates: $\{t,r, \vartheta \in \mathbb{S}^{n-1}\}, \quad g = -\ed t^2 + \ed r^2 + r^2  \circg$, where $\circg$ is the standard round metric on the sphere $\mathbb{S}^{n-1}$ of radius $1$;

\item double null coordinates: $\{u, \uu, \vartheta \in \mathbb{S}^{n-1}\}, \quad g = -4 \ed u \ed \uu + r^2  \circg$.

\end{enumerate}
The coordinates are related by the following transformations:
\begin{enumerate}[leftmargin=.5in]
\item[a.$\rightarrow$b.] 
$t=x^0, \quad r = \sqrt{(x^1)^2 + \cdots + (x^n)^2}, \quad \partial_t = \partial_0, \quad \partial_r = \sum_{i=1}^n \frac{x^i}{r}\partial_i$;
\item[b.$\rightarrow$c.] 
$u=\frac{t-r}{2}, \quad \uu = \frac{t+r}{2}, \quad \partial_u = \partial_t - \partial_r, \quad \partial_{\uu} = \partial_t + \partial_r$.
\end{enumerate}
In the double null coordinate system $\{u, \uu, \vartheta \in \mathbb{S}^{n-1}\}$, the level sets of $u$ and $v$ are lightcones in $\mathbb{M}^{n+1}$. We have the following 
\begin{enumerate}
\item[$C_u$:] 
Level set of $u$, which is the outgoing lightcone in $\mathbb{M}^{n+1}$;
\item[$\uC_{\uu}$:] 
Level set of $\uu$, which is the incoming lightcone in $\mathbb{M}^{n+1}$.
\end{enumerate}

\begin{figure}[H]
\centering
\begin{tikzpicture}[scale=0.8]
\draw[fill] (0,0) circle (1pt);
\draw (0,0) node[below right] {$o$}  -- (4,4) (0,0) -- (-4,4);
\draw[thick, dashed] (3,3) to [out=45,in=15] (2,3.5)
to [out=-165,in=-40] (-2,3)
to [out=140,in=135] (-3.5,3.5); 
\draw[thick] (-3.5,3.5) to [out=-45,in=-160] (-1.8,3.2)
to [out=20,in=-170] (-1,3.5)
to [out=10,in=150] (1,2.5)
to [out=-30,in=-135] (3,3); 
\draw (0,0) --(-1.8,3.2)
(0,0) -- (-1,3.5) node[above] {$f$} 
(0,0) -- (1,2.5);
\draw[dashed] (0,0) -- (2,3.5)  (0,0) -- (-2,3);
\draw[->] (2,3.6) to [out=60,in=180] (3,4) node[right] {$S$};
\end{tikzpicture}
\caption{Parameterisation of $S$}
\label{fig 2}
\end{figure}
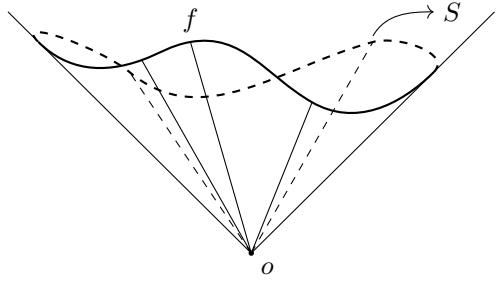

Considering the outgoing lightcone $C_{u=0}$ (simply denoted by $C_0$ from now on), $\{\uu \geq 0, \vartheta\in \mathbb{S}^{n-1}\}$ is a coordinate system of $C_0$. Then we can parameterise a spacelike surface $S$ in $C_0$ by a positive function $f$ as its graph of $\uu$ over the $\vartheta$ domain (see figure \ref{fig 2}):
\begin{align*}
S = \{ (v,\vartheta): v = f(\vartheta)>0\}.
\end{align*}
Let $\slashg$ be the intrinsic metric of $S$, which is $g|_{S}$, then $\slashg= f^2 \circg$.

\subsection{Geometry of a spacelike surface in a lightcone}\label{subsec 2.2}
A conjugate null frame $\{L, \uL\}$ of a spacelike surface $S$ consists of two future directed normal null vectors with the following condition
\begin{align*}
g( L, \uL )= -2.
\end{align*}
See figure \ref{fig  3}. Asscociated with the spacelike surface $S$ in $C_0$, we introduce the following special conjugate null frame $\{L, \uL\}$ in its normal bundle that
\begin{align*}
L= f \partial_{\uu},
\quad
g( L, \uL )= -2.
\end{align*}
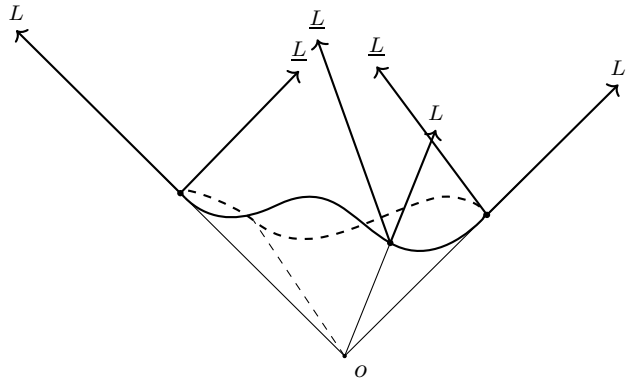
\begin{figure}[H]
\centering
\begin{tikzpicture}[scale=0.6]
\draw[fill] (0,0) circle (1pt);
\draw (0,0) node[below right] {$o$}  -- (4,4) (0,0) -- (-4,4);
\draw[thick, dashed] (3,3) to [out=45,in=15] (2,3.5)
to [out=-165,in=-40] (-2,3)
to [out=140,in=135] (-3.5,3.5); 
\draw[dashed] 
(0,0) -- (-2,3);
\draw[thick] (-3.5,3.5) to [out=-45,in=-160] (-1.8,3.2)
to [out=20,in=-170] (-1,3.5) 
to [out=10,in=150] (1,2.5)
to [out=-30,in=-135] (3,3); 
\draw 
(0,0) -- (1,2.5); 
\path[fill=black] (1,2.5) circle(2pt);
\draw[->,thick] (1,2.5) -- (2,5) node[above] {\footnotesize$L$};  
\draw[->,thick] (1,2.5) -- (-0.6,7) node[above] {\footnotesize$\uL$}; 
\path[fill=black] (3.12,3.12) circle(2pt);
\draw[->,thick] (3.12,3.12) -- (6,6) node[above] {\footnotesize$L$}; 
\draw[->,thick] (3.12,3.12) -- (0.7,6.4) node[above] {\footnotesize$\uL$}; 
\path[fill=black] (-3.6,3.6) circle(2pt);
\draw[->,thick] (-3.6,3.6) -- (-7.2,7.2) node[above] {\footnotesize$L$}; 
\draw[->,thick] (-3.6,3.6) -- (-1,6.3) node[above] {\footnotesize$\uL$}; 
\end{tikzpicture}
\caption{Conjugate null frame $\{ L, \uL \}$.}
\label{fig 3}
\end{figure}
Then we can define the second fundamental forms and the torsion relative to this conjugate null frame.
\begin{definition}\label{def 2.1}
Let $S$ be a spacelike surface in the outgoing lightcone $C_0$. Let $\{L, \uL\}$ be a conjugate null frame of $S$, then we define the second fundamental forms $\chi$, $\uchi$ and the torsion $\eta$ as follows: let $\nabla$ be the covariant derivative of $g$, then for $X,Y\in TS$,
\begin{align*}
\chi(X,Y) = g(\nabla_X L, Y),
\quad
\uchi(X,Y) = g(\nabla_X \uL, Y),
\quad
\eta(X) = \frac{1}{2} g(\nabla_X L, \uL).
\end{align*}
Decompose $\chi$ and $\uchi$ into trace part and trace-free part,
\begin{align*}
\chi = \hatchi + \frac{1}{n-1} \tr \chi \slashg,
\quad
\uchi = \hatuchi + \frac{1}{n-1} \tr \uchi \slashg.
\end{align*}
We call $\tr \chi$ the outgoing null expansion, $\tr \uchi$ the incoming null expansion, $\hatchi$ and $\hatuchi$ the outgoing and incoming shears respectively.
\end{definition}
We can easily calculate $\chi$ that $\chi= \slashg$ and $\tr \chi = n-1$ since $\nabla_X L = X$. Calculate $\uchi$ and $\eta$ in the following. Let $\{\theta^1, \cdots, \theta^{n-1}\}$ be a coordinate system of $\mathbb{S}^{n-1}$, then we introduce the following coordinate systems of $\mathbb{M}^{n+1}$ and $S$:
\begin{enumerate}[label=\alph*.]
\item 
$\{u, \uu, \theta^1, \cdots, \theta^{n-1}\}$: a double null coordinate system of $\mathbb{M}^{n+1}$, with the coordinate vectors 
\begin{align*}
\partial_u, \quad \partial_{\uu}, \quad \partial_i, \quad i=1,\cdots,n-1.
\end{align*}
The covariant derivatives of the coordinate vectors are given by
\begin{align*}
&
\nabla_{\partial_i} \partial_{\uu}= \frac{1}{r} \partial_i,
\quad
\nabla_{\partial_i} \partial_u = - \frac{1}{r} \partial_i,
\\
&
\nabla_{\partial_{\uu}} \partial_i =  \frac{1}{r} \partial_i,
\quad
\nabla_{\partial_u} \partial_i = - \frac{1}{r} \partial_i,
\\
&
\nabla_{\partial_u} \partial_{\uu} = \nabla_{\partial_{\uu}} \partial_u =0,
\quad
\nabla_{\partial_i} \partial_j = \circnabla_{\partial_i} \partial_j - \frac{r}{2} \circg_{ij} \partial_{\uu} + \frac{r}{2} \circg_{ij} \partial_u,
\end{align*}
where $\circnabla$ is the covariant derivative of $\circg$.

\item
$\{\theta^1, \cdots, \theta^{n-1}\}$: a coordinate system of $S$ when restricted on it, with the coordinate vectors
\begin{align*}
\dpartial_i = \partial_i + f_i \partial_{\uu},
\quad
f_i = \frac{\partial f}{\partial \theta^i}.
\end{align*}
\end{enumerate}
The conjugate null frame $\{L, \uL\}$ is given by
\begin{align*}
L = f \partial_{\uu},
\quad
\uL
= f^{-1} \big[ \partial_u + \big( f^{-2} \circg^{ij} f_i f_j \big) \partial_{\uu} + 2 \big( f^{-2} \circg^{ij} f_j \big) \partial_i \big].
\end{align*}
Thus we have the following results of $\uchi$ and $\eta$.
\begin{enumerate}
\item[$\uchi$:] $\uchi_{ij} = g(\nabla_{\dpartial_i} \uL, \dpartial_j) = - g(\nabla_{\dpartial_i} \dpartial_j, \uL)$. Since
\begin{align*}
\nabla_{\dpartial_i} \dpartial_j 
= 
\circnabla_{\partial_i} \partial_j - \frac{f}{2} \circg_{ij} \partial_{\uu} + \frac{f}{2} \circg_{ij} \partial_u
+ f_i f^{-1} \partial_j
+ f_{ij} \partial_{\uu} + f_j f^{-1} \partial_i,
\end{align*}
then
\begin{align}
\begin{aligned}
&
\uchi_{ij} = -g(\nabla_{\dpartial_i} \dpartial_j, \uL)
=
- \circg_{ij} + f^{-2} \circg^{kl} f_k f_l \circg_{ij} + 2f^{-1} \circnabla^2_{ij} f - 4 f^{-2} f_i f_j,
\\
&
\tr \uchi = -(n-1) f^{-2} + (n-5) f^{-4} \circg^{ij} f_i f_j + 2 f^{-3} \circDelta f.
\end{aligned}
\label{eqn 2.1}
\end{align}

\item[$\eta$:] Since $\nabla_{\dpartial_i} L = \dpartial_i$, thus
\begin{align*}
\eta_i = \frac{1}{2} g(\nabla_{\dpartial_i} L, \uL) =\frac{1}{2} g(\dpartial_i, \uL) =0.
\end{align*}
\end{enumerate}
Let $\R_{ijkl}$ be the curvature\footnote{The convention of the Riemann curvature operator here: $\mathrm{R}_{ijkl} = g(\nabla_{\partial_i} \nabla_{\partial_j} \partial_k - \nabla_{\partial_j} \nabla_{\partial_i} \partial_k ,\partial_l )$.} of $(S,\slashg)$. We have the Gauss equation
\begin{align}
\begin{aligned}
&
\begin{aligned}
\R_{ijkl} 
&= \frac{1}{2}( \chi_{ik} \uchi_{jl} + \uchi_{ik} \chi_{jl} - \chi_{jk} \uchi_{il} - \uchi_{jk} \chi_{il} )
\\
&
= \frac{1}{2}( \slashg_{ik} \uchi_{jl} + \uchi_{ik} \slashg_{jl} - \slashg_{jk} \uchi_{il} - \uchi_{jk} \slashg_{il} )
\end{aligned}
\\
\Rightarrow\quad&
\R_{ij} = - \frac{n-3}{2} \uchi_{ij} - \frac{1}{2} \tr\uchi \cdot \slashg_{ij},
\quad
\R = -(n-2) \tr \uchi.
\end{aligned}
\label{eqn 2.2}
\end{align}

\subsection{Functional $\calF$ and corresponding variational problem}\label{subsec 2.3}
In the following, we shall first introduce the functional $\calF$ of a closed spacelike surface in $C_0$, then give a geometric interpretation of $\calF$. We shall also obtain the formula of $\calF$ in terms of the parameterisation of $S$.

\begin{definition}\label{def 2.2}
Let $S$ be a closed spacelike surface in the lightcone $C_0$ in the Minkowski spacetime $\mathbb{M}^{n+1}$. Let $\{L ,\uL\}$ be a general conjugate null frame of $S$. We define the functional $\calF$ of $S$ by
\begin{align*}
\calF[S] = \int_S \tr \chi \cdot \tr \uchi\, \dvol_{\subslashg}.
\end{align*}
In particular for the conjugate null frame $\{L,\uL\}$ with $L = f \partial_{\uu}$, we have that
\begin{align*}
\calF[S] = \int_S (n-1) \tr \uchi\, \dvol_{\subslashg}.
\end{align*}
\end{definition}
From the Gauss equation \eqref{eqn 2.2}, we obtain the following geometric interpretation of the functional $\calF$.
\begin{enumerate}[label=\footnotesize\textbullet]
\item By the geometric meaning of $\tr \uchi$, $\int_S \tr \uchi\, \dvol_{\subslashg}$ measures the infinitesimal change of the area (length if $n=2$) of $S$ in the direction of $\uL$.

\item $n=3$, $\dim S=2$. We have that $\R = -\tr \uchi$. Let $K$ be the Gauss curvature of $(S,\slashg)$, then $\R = 2K$ and
\begin{align*}
\calF[S] =-\int_S 2 \R\, \dvol_{\subslashg} = -4 \int_S K\, \dvol_{\subslashg} = -8\pi.
\end{align*}

\item $n\geq 4$, $\dim S= n-1 \geq 3$. We have that $\R = -(n-2) \tr \uchi$, thus
\begin{align*}
\calF[S] = -\frac{n-1}{n-2} \int_S \R\, \dvol_{\subslashg}.
\end{align*}
Note that the integral $\int_S \R\, \dvol_{\subslashg}$ is the Einstein-Hilbert functional of $(S,\slashg)$.

\end{enumerate}

The functional $\calF[S]$ naturally reminds one of the Yamabe problem in the case of $\dim S =n-1 \geq 3$. We refer to the nice expository articles \cite{LP1987}\cite{S1989} for the thorough discussion of the Yamabe problem. Recall that the Yamabe problem asks that, for a Riemannian manifold $(M,g_0), \dim M \geq 3$, find the extreme point of the Einstein-Hilbert functional in the conformal class $\mathcal{M}_1$ of $g_0$ with the volume equal to $1$, i.e. let
\begin{align*}
&
I(M,g_0) = \inf\{ \int_M \R_{g}\, \dvol_g:  g \in \mathcal{M}_1 \},
\\
&
\mathcal{M}_1 = \{ g: g= e^{h} g_0, \mathrm{Vol}(M,g) = \int_M \dvol_g =1   \},
\end{align*}
find out $I(M,g_0)$ and whether that it can be achieved by a conformal metric $g\in \mathcal{M}_1$. Another equivalent formulation of the Yamabe problem is that find the conformal metric $g$ achieving the following Yamabe invariant $\lambda(M,g_0)$ of a Riemannian manifold $(M,g_0)$ with $\dim M =n-1 \geq 3$,
\begin{align*}
\lambda(M,g_0) = \inf\{ Q(M,g): g=e^h g_0\}
\quad
Q(M,g) = \frac{\int_M \R_g \dvol_g}{[\mathrm{Vol}(M,g)]^{\frac{n-3}{n-1}}}.
\end{align*}

Motived by the Yamabe problem, we can ask the following problem for the functional $\calF$.
\begin{problem}\label{pro 2.3}
In the Minkowski spacetime $(\mathbb{M}^{n+1},g), n\geq 2$, find the extreme point of the following functional
\begin{align*}
\mathcal{Q}(S) = \frac{\calF[S]}{|S|^{\frac{n-3}{n-1}}} = \frac{ \int_S \tr \chi \cdot \tr \uchi\, \dvol_g}{\big[\int_S \dvol_{\subslashg} \big]^{\frac{n-3}{n-1}}}.
\end{align*}
\end{problem}

We already have the answer of the above problem in the case $n\geq 3$ as follows.
\begin{enumerate}[label=\footnotesize\textbullet]
\item $n=3, \dim S=2$. $\mathcal{Q}(S) = -8\pi$ is a constant. The problem is trivial.

\item $n\geq 4, \dim S =n-1 \geq 3$. The problem is equivalent to the Yamabe problem of the standard round sphere $(\mathbb{S}^{n-1} ,\circg)$. We have that 
\begin{align*}
&\textstyle
\sup\{ \mathcal{Q}[S] \} 
=
-\frac{n-1}{n-2} \inf\{ Q(\mathbb{S}^{n-1}, e^h \circg) \}
=
-\frac{n-1}{n-2} \lambda (\mathbb{S}^{n-1}, \circg)
=
- (n-1)^2\, \big[ \mathrm{Vol}(\mathbb{S}^{n-1} ,\circg) \big]^{\frac{2}{n-1}},
\\
&\textstyle
\mathrm{Vol}(\mathbb{S}^{n-1} ,\circg) 
= 
\frac{2\pi^{\frac{n}{2}}}{\Gamma(\frac{n}{2})},
\end{align*}
and $\sup\{ \mathcal{Q} [S] \}$ is achieved by the spacelike hyperplane section of the lightcone $C_0$.
\end{enumerate}

The remaining case of the problem is when $n=2$, $S$ is a spacelike curve in the lightcone $C_0$. The problem is reduced to a variational problem on the circle. We derive an explicit formula of this problem in the following.
\begin{enumerate}[label=\roman*.]
\item Let $\theta$ be the parameter of the circle. Any closed spacelike curve $S$ in $C_0$ can be parametrised by a positive function $f$ on the circle. Then the metric and volume form on $S$ are given by
\begin{align*}
\slashg = f^2 \ed \theta^2,
\quad 
\dvol_{\subslashg} = f \ed \theta.
\end{align*}

\item By equation \eqref{eqn 2.1}, we have
\begin{align*}
\tr \uchi
= - f^{-2}  - 3 f^{-4} (f_{\theta})^2 + 2 f^{-3} f_{\theta \theta}.
\end{align*}

\item Abuse the notation that $\calF[f] = \calF[S]$ and $\mathcal{Q}[f] = \mathcal{Q}[S]$. Then
\begin{align}
\begin{aligned}
\calF[f] 
&
= \int_S \tr \uchi\, \dvol_{\subslashg} 
= \int_0^{2 \pi} \big[ -f^{-1} + f^{-3}(f_{\theta})^2 \big] \ed \theta,
\\
\calQ[f]
&
= \calF[f] \cdot |S|
= \int_0^{2 \pi} \big[ -f^{-1} + f^{-3}(f_{\theta})^2 \big]\, \ed \theta \cdot \int_0^{2 \pi} f \ed \theta.
\end{aligned}
\label{eqn 2.3}
\end{align}

\item We introduce two different functions by $h = \log f$ and $v = f^{-\frac{1}{2}} = e^{-\frac{h}{2}}$, and abuse the notations that $\calF[f] = \calF[h] = \calF[v]$ and $\calQ[f] = \calQ[h] = \calQ[v]$, then
\begin{align}
&
\calF[h]
= \int_0^{2 \pi} e^{-h} \big[  (h_{\theta})^2 -1 \big] \ed \theta,
&&
\calQ[h]
= \int_0^{2 \pi} e^{-h} \big[ (h_{\theta})^2 -1 \big] \ed \theta \cdot \int_0^{2 \pi} e^{h} \ed \theta,
\label{eqn 2.4}
\\
&
\calF[v]
= \int_0^{2 \pi} \big[ 4(v_{\theta})^2 - v^2 \big] \ed \theta,
&&
\calQ[v]
= \int_0^{2 \pi} \big[ 4(v_{\theta})^2 - v^2 \big] \ed \theta \cdot \int_0^{2 \pi} v^{-2} \ed \theta,
\label{eqn 2.5}
\end{align}

\end{enumerate}

We can reformulate problem \ref{pro 2.3} as follows.
\begin{problem}\label{pro 2.4}
In the Minkowski spacetime $(\mathbb{M}^{3}, g)$, find the extreme point of $\calQ[f]=\calQ[h] = \calQ[v]$, where $f= e^h = v^{-2}$ is the parameterisation of a spacelike curve in a lightcone $C_0$.

The above variational problem of $\calQ$ is equivalent to the variational problem of finding the extreme point of $\calF$ under the constraint $|S| = \int_0^{2\pi} f \ed \theta = \int_0^{2\pi} e^h \ed \theta = \int_0^{2\pi} v^{-2} \ed \theta = 2 \pi$.
\end{problem}

\section{Lorentz invariance of functional $\calF$}\label{sec 3}
From the geometric meaning of $\calF$ and $\calQ$ in definition \ref{def 2.2} and problem \ref{pro 2.3}, clearly that they are invariant under the Lorentzian group of the Minkowski spacetime $(\mathbb{M}^3, g)$, i.e. let $S$ be an arbitrary spacelike curve in $C_0$, then we have that
\begin{align*}
\calF[\varphi(S)] = \calF[S],
\quad
\calQ[\varphi(S)] = \calQ[S],
\quad
\text{for every Lorentz transformation $\varphi$ of $(\mathbb{M}^3, g)$.}
\end{align*}

In the following, we work out the explicit formula of the Lorentz transformation on the parameterisation function $f$ of $S$. The Lorentzian group of $(\mathbb{M}^3,g)$ is $3$ dimensional. It can be parameterised by three parameters $\alpha \in (-1,1), \theta_0, \btheta_0 \in \mathbb{R}$ that in the spatial polar coordinate system $(t, r, \theta)$, the corresponding Lorentz transformation $\varphi_{\alpha, \theta_0, \btheta_0}$ takes the form
\begin{align*}
\varphi_{\alpha, \theta_0, \btheta_0}: 
\quad
(t,r,\theta)  
\quad
\mapsto 
\quad
(\bt, \br, \btheta),
\end{align*}
where
\begin{align*}
\bt = \frac{t + \alpha r \cos(\theta-\theta_0)}{\sqrt{1- \alpha^2}},
\quad
\br =\sqrt{\frac{\big(r+  \alpha t \cos(\theta-\theta_0)\big)^2 - \alpha^2 (r^2-t^2) \sin^2 (\theta-\theta_0)}{1- \alpha^2}},
\end{align*}
and $\btheta$ is solved from the following equations
\begin{align*}
\cos (\btheta - \btheta_0)
=
\frac{\alpha t+ r\cos (\theta-\theta_0)}{\sqrt{\big(r+ \alpha t \cos(\theta-\theta_0)\big)^2 - \alpha^2 (r^2-t^2) \sin^2 (\theta-\theta_0)}},
\\
\sin (\btheta - \btheta_0)
=
\frac{r \sin(\theta-\theta_0) \sqrt{1-\alpha^2} }{ \sqrt{\big(r+ \alpha t \cos(\theta-\theta_0)\big)^2 - \alpha^2 (r^2-t^2) \sin^2 (\theta-\theta_0)} }.
\end{align*} 
The restriction of $\varphi_{\alpha, \theta_0, \btheta_0}$ on $C_0$ in the coordinate system $(\uu, \theta)$ is
\begin{align*}
\varphi_{\alpha, \theta_0, \btheta_0}:
\quad
&
(\uu, \theta)
\quad
\mapsto
\quad
(\buu, \btheta)
= 
\Big( \frac{1+\alpha \cos(\theta-\theta_0)}{\sqrt{1-\alpha^2}} \uu, \btheta \Big),
\\
\varphi_{-\alpha, \btheta_0, \theta_0}
=
\varphi_{\alpha, \theta_0, \btheta_0}^{-1}:
\quad
&
(\buu, \btheta)
\quad
\mapsto
\quad
(\uu, \theta)
= 
\Big( \frac{1-\alpha \cos (\btheta-\btheta_0)}{\sqrt{1-\alpha^2}} \buu, \theta \Big),
\end{align*}
where $\btheta$ is solved from the equations
\begin{align}
\begin{aligned}
&
\cos  (\btheta - \btheta_0)
=
\frac{\alpha + \cos (\theta-\theta_0)}{1+ \alpha  \cos(\theta-\theta_0)},
\quad
\sin  (\btheta - \btheta_0)
=
\frac{  \sqrt{1-\alpha^2} \sin(\theta-\theta_0)}{ 1+ \alpha  \cos(\theta-\theta_0) },
\\
\Leftrightarrow
\quad
&
\cos (\theta-\theta_0) = \frac{\cos  (\btheta - \btheta_0) - \alpha}{1 - \alpha \cos  (\btheta - \btheta_0)},
\quad
\sin (\theta-\theta_0) = \frac{ \sqrt{1-\alpha^2} \sin  (\btheta - \btheta_0) }{1 - \alpha \cos (\btheta - \btheta_0)}.
\end{aligned}
\label{eqn 3.1}
\end{align}
Let $\barf$ be the parameterisation function of $\varphi_{\alpha, \theta_o, \btheta_0} (S)$, then we have the following formula of $\barf$,
\begin{align*}
\barf(\btheta) 
= \frac{1+\alpha \cos(\theta-\theta_0)}{\sqrt{1-\alpha^2}} f(\theta) 
=
\frac{\sqrt{1-\alpha^2}}{1- \alpha \cos (\btheta - \btheta_0)} f(\theta),
\end{align*}
where $\theta$ is determined by equations \eqref{eqn 3.1}. We introduce the notation $\nu_{\alpha, \btheta_0}$ to denoted the multiplier function
\begin{align*}
\nu_{\alpha, \btheta_0} (\btheta) = \frac{\sqrt{1-\alpha^2}}{1- \alpha \cos (\btheta - \btheta_0)} .
\end{align*}
We can state the Lorentz invariance of the functionals $\calF$ and $\calQ$ now.
\begin{proposition}\label{pro 3.1}
Let $S$ be a spacelike curve in $C_0$ of $(\mathbb{M}^3,g)$. Let $f= e^{h} = v^{-2}$ be the parameterisation function of $S$. Let $\varphi_{\alpha,\theta_0,\btheta_0}$ be the Lortenzian transformation. Then
\begin{align*}
\calF[\varphi(S)] = \calF[S],
\quad
\calQ[\varphi(S)] = \calQ[S].
\end{align*}
Let $\barf = e^{\bh} = \bv^{-2}$ be the parameterisation function of $\varphi_{\alpha, \theta_0, \btheta_0}(S)$. Then $\varphi_{\alpha,\theta_0, \btheta_0}$ induces a transformation from the parameterisation functions $f$ to $\barf = \varphi_{\alpha,\theta_0, \btheta_0}(f)$,
\begin{align*}
\barf(\btheta) = \nu_{\alpha, \btheta_0}(\btheta) f(\theta),
\quad
\text{where $\theta$ and $\btheta$ are related by equations \eqref{eqn 3.1}.}
\end{align*}
Then we have
\begin{align*}
\calF[\barf]=\calF[\varphi_{\alpha,\theta_0, \btheta_0}(f)] = \calF[f],
\quad
\calQ[\barf]=\calQ[\varphi_{\alpha,\theta_0, \btheta_0}(f)] = \calQ[f].
\end{align*}
In terms of the parameterisation functions $h, v$ and $\bh = \rho_{\alpha,\theta_0, \btheta_0}(h), \bv=\gamma_{\alpha,\theta_0, \btheta_0}(v)$, where the transformation takes the form
\begin{align*}
\bh(\btheta) = h(\theta) + \log \nu_{\alpha, \btheta_0}(\btheta),
\quad
\bv(\btheta) = [\nu_{\alpha, \btheta_0}(\btheta)]^{-\frac{1}{2}} v(\theta).
\end{align*}
Then we have
\begin{align*}
&
\calF[\bh]=\calF[\rho_{\alpha,\theta_0, \btheta_0}(h)]= \calF[h],
&&
\calQ[\bh]=\calQ[\rho_{\alpha,\theta_0, \btheta_0}(h)]= \calQ[h],
\\
&
\calF[\bv]=\calF[\gamma_{\alpha,\theta_0, \btheta_0}(v)] = \calF[v],
&&
\calQ[\bv]=\calQ[\gamma_{\alpha,\theta_0, \btheta_0}(v)] = \calQ[v].
\end{align*}
\end{proposition}

\section{Critical points of $\calF$ under constraint and local analysis at critical points}\label{sec 4}
In this section, we derive the Euler-Lagrangian equation for the critical points of $\calF$ under the constraint $|S|=2\pi$ and solve the equation. We calculate the second variation of $\calF[S]$ at the critical points and deduce that they are local minimisers of $\calF[S]$ under the constraint.

\subsection{Euler-Lagrangian equation of $\calF$ under constraint}\label{subsec 4.1}
Let $S_t$ be a family of spacelike curves in $C_0$ satisfies the constraint $|S_t|=2\pi$. Let $f_t$ be the parameterisation function of $S_t$, then
\begin{align*}
\int_0^{2\pi} f_t\, \ed \theta=2\pi.
\end{align*}
In the following, we shall use the dot $\cdot$ to denote the derivative with respect to $t$ and the prime $'$ to denote the derivative with respect to $\theta$. Then
\begin{align*}
\int_0^{2\pi} \dot{f_t}\, \ed \theta =0.
\end{align*}
We calculate the first variation of $\calF[f_t]$.
\begin{align*}
\frac{\ed}{\ed t} \calF[f_t]
&
= 
\frac{\ed}{\ed t} \int_0^{2 \pi} \big[ -f_t^{-1} + f_t^{-3}(f_t')^2 \big] \ed \theta
\\
&
=
\int_0^{2 \pi} \big[ f_t^{-2} \dot{f_t} -3 f_t^{-4}(f_t')^2 \dot{f_t} + 2 f_t^{-3} f_t' \dot{f_t}' \big] \ed \theta
\\
&
=
\int_0^{2 \pi} \big[ f_t^{-2} + 3 f_t^{-4}(f_t')^2 - 2 f_t^{-3} f_t'' \big] \dot{f_t}  \ed \theta.
\end{align*}
Thus the critical point of $\calF$ under the constraint satisfies the Euler-Lagrangian equation
\begin{align}
f^{-2} + 3 f^{-4}(f')^2 - 2 f^{-3} f'' = \lambda,
\quad
\lambda \in \mathbb{R}.
\label{eqn 4.1}
\end{align}

\subsection{Critical points of $\calF$ under constraint}\label{subsec 4.2}
We solve the above equation in the following. From equation \eqref{eqn 4.1},
\begin{align*}
[f^{-3} (f')^2]'  + (f^{-1})'  + \lambda f' =0
\quad
\Rightarrow
\quad
(f')^2 = -f^2 + \mu f^3 - \lambda f^4.
\end{align*}
Thus
\begin{align*}
\int \frac{\ed f}{f\sqrt{-\lambda f^2 + \mu f -1}} = \theta.
\end{align*}
Since we want a positive solution of $f$ on the circle, the quadratic equation $-\lambda f^2 + \mu f -1$ shall have two real roots (multiplicity counted). If two real roots are equal, then $f$ is a constant function on the circle. In the following, we consider the case of two distinct roots. Suppose that
\begin{align*}
f_{\min} = \frac{\mu - \sqrt{\mu^2-4\lambda}}{2\lambda} >0,
\quad
f_{\max} = \frac{\mu + \sqrt{\mu^2-4\lambda}}{2\lambda},
\quad
\mu^2-4\lambda>0,
\,
\mu>0,
\,
\lambda>0,
\end{align*}
In order to solve $f$, we introduce the notations
\begin{align*}
&
\alpha = \sqrt{\frac{\mu^2 - 4\lambda}{\mu^2}} < 1,
\quad
f_{\min} = \frac{1}{\sqrt{\lambda}} \sqrt{\frac{1-\alpha}{1+\alpha}},
\quad
f_{\max} = \frac{1}{\sqrt{\lambda}} \sqrt{\frac{1+\alpha}{1-\alpha}},
\\
&
\btheta = \pi - \arccos \Big( \frac{\sqrt{\lambda} \sqrt{1-\alpha^2} f -1}{\alpha} \Big),
\end{align*}
then
\begin{align*}
\int \frac{\sqrt{1-\alpha^2}}{1 - \alpha \cos \btheta} \ed \btheta = \theta
\quad
\Rightarrow
\quad
\tan \frac{\btheta}{2} = \sqrt{\frac{1-\alpha}{1+\alpha}} \tan \frac{\theta - \theta_0}{2}.
\end{align*}
Hence
\begin{align*}
f = \frac{1}{\sqrt{\lambda}} \frac{\sqrt{1-\alpha^2}}{1+ \alpha \cos(\theta - \theta_0)}.
\end{align*}
We determine the value of $\lambda$ to be $1$ from the constraint
\begin{align*}
2\pi = \int_0^{2\pi} f\, \ed \theta = \int_0^{2\pi} \frac{1}{\sqrt{\lambda}} \frac{\sqrt{1-\alpha^2}}{1+ \alpha \cos(\theta - \theta_0)} \ed \theta = \frac{2\pi}{\sqrt{\lambda}}.
\end{align*}
\begin{proposition}\label{pro 4.1}
The set of critical points of the functional $\calF$ under the constraint $\int_0^{2\pi} f \ed \theta =2\pi$ consists of the following functions
\begin{align*}
f(\theta) = e^{h(\theta)} = v^{-2}(\theta)
=
\frac{\sqrt{1-\alpha^2}}{1+ \alpha \cos(\theta - \theta_0)}
=
\nu_{-\alpha,\theta_0}(\theta),
\quad
\alpha \in (-1,1),
\quad
\theta_0 \in \mathbb{R}.
\end{align*}
These critical points are all related by the Lorentz transformation in proposition \ref{pro 3.1}.
\end{proposition}

\subsection{Second variation of $\calF$ at critical points under constraint}\label{subsec 4.3}
In the following, we shall mainly work with the parameterisation function $v$. We calculate the second variation of the functional $\calF[v] = \int_0^{2\pi} [4 (v')^2 - v^2] \ed \theta$ at $v_0=1$ under the constraint $\int_0^{2\pi} v^{-2} \ed \theta=2\pi$. Let $v_t$ be a variation of $v_0$ under the constraint. We have that
\begin{align*}
\frac{\ed}{\ed t} \int_0^{2\pi} v_t^{-2} \ed \theta
&=
\int_0^{2\pi} -2 \dot{v}_t \cdot v_t^{-3} \ed \theta 
= 
0,
\\
\frac{\ed^2}{\ed t^2} \int_0^{2\pi} v_t^{-2} \ed \theta
&=
\int_0^{2\pi} -2 (\ddot{v}_t v_t -3 \dot{v}_t^2 ) \cdot v_t^{-4} \ed \theta 
= 
0,
\\
\frac{\ed}{\ed t} \calF[v_t] 
&= 
- 2 \int_0^{2\pi} (4 v_t'' + v_t ) \cdot \dot{v}_t \, \ed \theta,
\\
\frac{\ed^2}{\ed t^2} \calF[v_t] 
&= 
- 2 \int_0^{2\pi} (4 v_t'' + v_t ) \cdot \ddot{v}_t \, \ed \theta
- 2 \int_0^{2\pi} (4 \dot{v}_t'' + \dot{v}_t ) \cdot \dot{v}_t \, \ed \theta,
\end{align*}
Thus at $t=0$ we have
\begin{align*}
\frac{\ed}{\ed t} \Big|_{t=0} \int_0^{2\pi} v_t^{-2} \ed \theta
&=
\int_0^{2\pi} -2 \dot{v}\, \ed \theta 
= 
0,
\\
\frac{\ed^2}{\ed t^2} \Big|_{t=0} \int_0^{2\pi} v_t^{-2} \ed \theta
&=
\int_0^{2\pi} -2 (\ddot{v} -3 \dot{v}^2 ) \ed \theta 
= 
0,
\\
\frac{\ed}{\ed t} \Big|_{t=0} \calF[v_t]
&= 
- 2 \int_0^{2\pi} \dot{v} \, \ed \theta
=0,
\\
\frac{\ed^2}{\ed t^2} \Big|_{t=0} \calF[v_t]
&= 
- 2 \int_0^{2\pi} \ddot{v} \, \ed \theta
- 2 \int_0^{2\pi} (4 \dot{v}'' + \dot{v}) \cdot \dot{v} \, \ed \theta
=
8 \int_0^{2\pi} [ (\dot{v}')^2 - \dot{v}^2 ] \ed \theta.
\end{align*}
Thus we conclude that the second variation of $\calF[v]$ at $v_0=1$ under the constraint $\int_0^{2\pi} v^{-2} \ed \theta= 2\pi$ is nonnegative with the null space $\mathrm{span}\{ \sin \theta, \cos \theta\}$. Define the operator $\calL$
\begin{align*}
\calL v = -8v'' - 8v,
\end{align*}
then the second variation of $\calF$ takes the form
\begin{align*}
\frac{\ed^2}{\ed t^2} \Big|_{t=0} \calF[v_t] = \int_0^{2\pi} \dot{v} \cdot \calL \dot{v} \ed \theta.
\end{align*}
The eigenvalues and eigenfunctions of $\calL$ in the space $\{ v: \int_0^{2\pi} v\, \ed \theta=0 \}$ are
\begin{align*}
\kappa_m = 8m(m+2),
\quad
V_m = \mathrm{span}\{ \sin \big( (m+1) \theta \big), \cos \big( (m+1) \theta \big) \},
\quad
n\in \mathbb{N}
\end{align*}
Note that the null space $V_0$ corresponds to the variation of $v_0=1$ through the critical points $v^{-2} = \nu_{\alpha,\theta_0}$. From the above local analysis of $\calF[v]$ at $v_0=1$, we conclude the following proposition.
\begin{proposition}\label{pro 4.2}
The critical points $v^{-2} = \nu_{\alpha, \theta_0}$ are local minimisers of $\calF[v]=\int_0^{2\pi} [4 (v')^2 - v^2] \ed \theta$ under the constraint $\int_0^{2\pi} v^{-2} \ed \theta =2\pi$. The value of the local minimum is $-2\pi$.

Equivalently, the critical points $v^{-2} = a \nu_{\alpha, \theta_0}, a>0$ are local minimisers of $\calQ[v] = \int_0^{2 \pi} \big[ 4(v')^2 - v^2 \big] \ed \theta \cdot \int_0^{2 \pi} v^{-2} \ed \theta$. The value of the local minimum is $-4\pi^2$.
\end{proposition}

From propositions \ref{pro 4.1} and \ref{pro 4.2}, we see that all the critical points are the same up to Lorentz transformations, and they are local minimisers. Thus it is natural to conjecture the following claim.
\begin{claim}\label{cla 4.3}
The critical points $v^{-2} = \nu_{\alpha, \theta_0}$ are actually global minimisers of $\calF[v]=\int_0^{2\pi} [4 (v')^2 - v^2] \ed \theta$ under the constraint $\int_0^{2\pi} v^{-2} \ed \theta =2\pi$, which is equivalent to the following inequality:
\begin{align*}
\int_{\mathbb{S}^1} [4(v')^2 - v^2]  \ed \theta \geq  - \frac{4\pi^2}{\int_{\mathbb{S}^1} v^{-2} \ed \theta}.
\end{align*}
\end{claim}
One possible way to prove this claim is to apply the min-max method and the mountain pass theorem: suppose that there exists a point $v_1$ satisfying the constraint $\int_0^{2\pi} v^{-2} \ed \theta =2\pi$ such that $\calF[v]<-2\pi$, then one could find the critical point by the min-max method to the set of curves from $v_0=1$ to $v_1$,
\begin{align*}
\min_{\gamma} \max_{t\in[0,1]} \Big\{\calF[\gamma(t)]: \gamma(0) = v_0, \gamma(1) = v_1, \gamma(t) = v_t, \int_0^{2\pi} v_t^{-2} \ed \theta=2\pi \Big\}.
\end{align*}
However there is a difficulty in the above approach, which is that the functional $\calF$ doesnot satisfy the Palais-Smale compactness condition, due to the fact that the invariance group of $\calF$ isnot compact. Thus there needs more work to answer whether $v^{-2}=\nu_{\alpha,\theta}$ are global minimisers of $\calF$ under the constraint using this approach.

In this paper, we explore another method to prove claim \ref{cla 4.3}, which employs the Lorentz transformation and the symmetric decreasing rearrangement. This approach is similar to the method of competing symmetries in the proof of the sharp Hardy-Littlewood-Sobolev inequality in \cite{CL1990}, where we also find a way to pick up one special minimiser $\nu_{\alpha,\theta_0}$ to break the symmetry of Lorentz transformation.

\section{Symmetric decreasing rearrangement and rearrangement inequality}\label{sec 5}
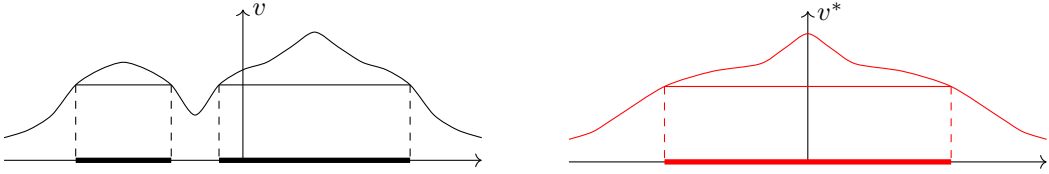
\begin{figure}[H]
\begin{subfigure}[h]{0.5\textwidth}
\centering
\begin{tikzpicture}
\draw[->] (0,0) -- (2*pi,0);
\draw[->] (1*pi,0) -- (1*pi,2) node[right] {$v$};
\draw plot [smooth] coordinates 
{(0,0+0.3) (0.1*pi,0.1+0.3) (0.2*pi,0.3+0.3) (0.3*pi,0.7+0.3) (0.4*pi,0.9+0.3) (0.5*pi,1+0.3) (0.6*pi,0.9+0.3) (0.7*pi,0.7+0.3) (0.8*pi,0.3+0.3) (0.9*pi,0.7+0.3) (1*pi,0.9+0.3) (1.1*pi,1+0.3) (1.2*pi,1.2+0.3) (1.3*pi,1.4+0.3) (1.4*pi,1.2+0.3) (1.5*pi,1+0.3) (1.6*pi,0.9+0.3) (1.7*pi,0.7+0.3) (1.8*pi,0.3+0.3) (1.9*pi,0.1+0.3) (2*pi,0+0.3)};
\draw (0.3*pi,0.7+0.3) -- (0.7*pi,0.7+0.3) (0.9*pi,0.7+0.3) -- (1.7*pi,0.7+0.3);
\draw[dashed] (0.3*pi,0.7+0.3) -- (0.3*pi,0) (0.7*pi,0.7+0.3) -- (0.7*pi,0) (0.9*pi,0.7+0.3) -- (0.9*pi,0) (1.7*pi,0.7+0.3) -- (1.7*pi,0);
\draw[line width=2pt] (0.3*pi,0) -- (0.7*pi,0) (0.9*pi,0) -- (1.7*pi,0);
\end{tikzpicture}
\end{subfigure}
\begin{subfigure}[h]{0.5\textwidth}
\centering
\begin{tikzpicture}
\draw[->] (0,0) -- (2*pi,0);
\draw[->] (1*pi,0) -- (1*pi,2) node[right] {$v^*$};
\draw[red] plot [smooth] coordinates 
{(0,0+0.3) (0.1*pi,0.1+0.3) (0.2*pi,0.3+0.3) (0.4*pi,0.7+0.3) (0.6*pi,0.9+0.3) (0.8*pi,1+0.3) (0.9*pi,1.2+0.3) (1*pi,1.4+0.3) (1.1*pi,1.2+0.3) (1.2*pi,1+0.3) (1.4*pi,0.9+0.3) (1.6*pi,0.7+0.3) (1.8*pi,0.3+0.3) (1.9*pi,0.1+0.3) (2*pi,0+0.3)};
\draw[red] (0.4*pi,0.7+0.3) -- (1.6*pi,0.7+0.3);
\draw[red] (0.4*pi,0.7+0.3) -- (1.6*pi,0.7+0.3);
\draw[dashed,red] (0.4*pi,0.7+0.3) -- (0.4*pi,0) (1.6*pi,0.7+0.3) -- (1.6*pi,0);
\draw[line width=2pt,red] (0.4*pi,0) -- (1.6*pi,0);
\end{tikzpicture}
\end{subfigure}
\caption{symmetric decreasing rearrangement}
\label{fig 4}
\end{figure}

In this section, we review the rearrangement inequality on the circle. Let $v$ be a function on the circle. Introduce the following notation
\begin{align*}
L_v(t) = \{\theta\in [-\pi,\pi]: v(\theta) > t\}
\end{align*}
We define the symmetric decreasing rearrangement $v^*$ of $v$ by symmetrising the level set of $v$ as follows. $v^*$ is monotone nonincreasing on $[0, \pi]$ such that
\begin{align*}
v^*(\theta) = v^*(-\theta),
\quad
| L_v(t) | = | L_{v^*} (t)|.
\end{align*}
See figure \ref{fig 4}. In another equivalent formulation,
\begin{align*}
v^*(\theta) = \sup \{t: |L_v(t)| > 2| \theta|  \}.
\end{align*}
$v^*$ is lower semicontinuous. If $v$ is continuous, then its symmetric decreasing rearrangement $v^*$ is also continuous. Moreover, $v^*$ is uniformly continuous with the same modulus of continuity as $v$.

We have the following identities and inequality for $v$ and $v^*$:
\begin{align*}
&
\int_0^{2\pi} v^2 \ed \theta = \int_0^{2\pi} (v^*\makebox[0ex]{})^2 \ed \theta,
\quad
\int_0^{2\pi} v^{-2} \ed \theta = \int_0^{2\pi} (v^*\makebox[0ex]{})^{-2} \ed \theta,
\\
&
\int_0^{2\pi} (v')^2 \ed \theta \geq \int_0^{2\pi} (v^*\makebox[0ex]{}')^2 \ed \theta.
\end{align*}
Then we have the following lemma on the symmetric decreasing rearrangement and the functional $\calF$.
\begin{lemma}\label{lem 5.1}
The symmetric decreasing rearrangement doesnot increase the functional $\calF$ under the constraint, i.e. for a positive function $v$ on the circle and its symmetric decreasing rearrangement $v^*$, we have that
\begin{align*}
\calF[v^*] \leq \calF[v],
\quad
\int_0^{2\pi} (v^*\makebox[0ex]{})^{-2} \ed \theta = \int_0^{2\pi} v^{-2} \ed \theta.
\end{align*}
\end{lemma}

\section{Decreasing sequence of $\calF$ with nondecreasing minimum}\label{sec 6}
In this section, we introduce a method to construct a decreasing sequence of the functional $\calF$ with nondecreasing minimum, using the symmetric decreasing rearrangement and the Lorentz transformation.

\subsection{Increase minimum of function by Lorentz transformation}\label{subsec 6.1}
In this subsection, we prove two lemmas on the minimum/infimum of a function and the Lorentz transformation, as preparations to construct the decreasing sequence of $\calF$. The lemmas are proved for functions on the circle and in the $3$-dim Minkowski spacetime here, but they can be easily generalised to the case of functions on $\mathbb{S}^{n-1}$ and the $(n+1)$-dim Minkowski spacetime.

For the sake of generality, most results in this section are proved for positive lower semicontinuous functions which makes the proof a bit more involved. One can simply assume that all the functions in the assumptions of the lemmas are continuous when going through for the first time.

\begin{lemma}\label{lem 6.1}
Let $v$ be a positive lower semicontinuous function on the circle satisfying $v(\theta) = v(-\theta)$ and monotonically nonincreasing on $[0,\pi]$, i.e. $v=v^*$. If $v^{-2}(\pi)> \lim_{x \searrow \frac{\pi}{2}} v^{-2}(x) = v^{-2} (\frac{\pi}{2})$, then there exists $\alpha_0 \in (0,1)$ such that for every Lorentz transformation $\varphi_{\alpha, 0, 0}, \alpha\in (0,\alpha_0]$ we have that \footnote{Note that $v$ is monotonically nonincreasing in $[0,\pi]$ and positive lower semicontinuous implies that $v^{-2}$ is monotonically nondecreasing in $[0,\pi]$ and positive upper semicontinuous. Thus $\varlimsup_{x \rightarrow \theta} v^{-2}(x) = \lim_{x \searrow \theta} v^{-2}(x)$, $v^{-2} (\theta) = \lim_{x \searrow \theta} v^{-2}(x) = v^{-2} (\theta^+)$ for $\theta \in [0,\pi]$ and $\max_{\theta} \{\varphi_{\alpha,0,0}(v^{-2})(\theta)\}$ exists. }
\begin{align*}
\max_{\theta} \{ \varphi_{\alpha, 0, 0} (v^{-2}) (\theta) \} < \max_{\theta} \{ v^{-2} (\theta) \},
\quad
\min_{\theta} \{ \gamma_{\alpha, 0, 0} (v) (\theta) \} > \min_{\theta} \{ v (\theta) \}.
\end{align*}
\end{lemma}
\begin{proof}
The lemma is easily observed from the geometric point of view, as illustrated in figure \ref{fig 5}. For $\alpha$ sufficiently small, the Lorentz transformation $\varphi_{\alpha,0,0}$ decrease the function $v^{-2}$ in the interval $[\frac{\pi}{2}+\epsilon, \pi]$ for some small $\epsilon$ (depending on $\alpha$), thus we can apply $\varphi_{\alpha,0,0}$ to decrease the maximum of $v$.
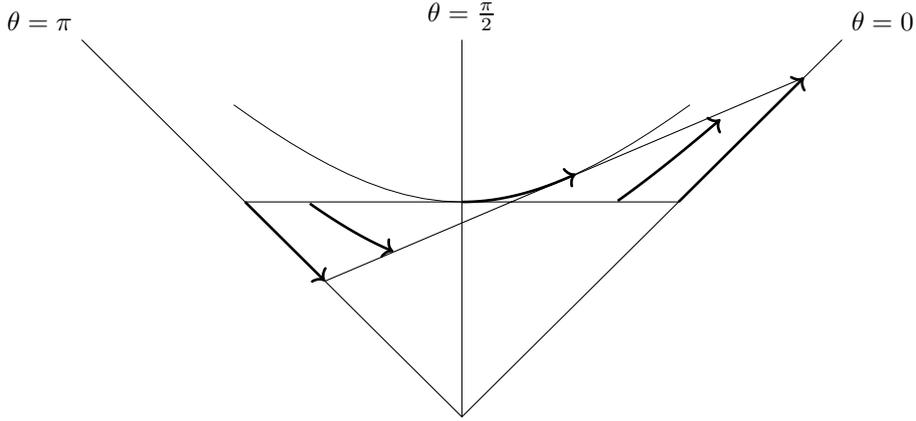
\begin{figure}[H]
\centering
\begin{tikzpicture}[scale=1]
\draw (-5,5) node[above left] {$\theta=\pi$} -- (0,0) -- (5,5) node[above right] {$\theta=0$};
\draw (0,0) -- (0,5) node[above] {$\theta =\frac{\pi}{2}$};
\draw [line width=0.2pt,domain=-3:3] plot (\x, {sqrt((2.85)^2+(\x)^2)});
\draw [<-,line width=1pt,domain=1.5:0] plot (\x, {sqrt((2.85)^2+(\x)^2)});
\draw [<-,line width=1pt,domain=-0.9:-2] plot (\x, {sqrt((2)^2+(\x)^2)});
\draw [<-,line width=1pt,domain=3.4:2.05] plot (\x, {sqrt((2)^2+(\x)^2)});
\draw (4.5,4.5) -- (-1.8,1.8) (-2.85,2.85) -- (2.85,2.85);
\draw[->,line width=1pt] (-2.85,2.85) -- (-1.8,1.8);
\draw[->,line width=1pt] (2.85,2.85) -- (4.5,4.5);
\end{tikzpicture}
\caption{Lorentzian transformation $\varphi_{\alpha,0,0}$ on $v^{-2}$.}
\label{fig 5}
\end{figure}
We prove the lemma by the above idea rigorously. Recall the formula of $\varphi_{\alpha,0,0}(v^{-2})$ that
\begin{align*}
\varphi_{\alpha,0,0}(v^{-2}) (\btheta) = \frac{\sqrt{1-\alpha^2}}{1 - \alpha \cos \btheta} v^{-2}(\theta),
\quad
\cos \theta = \frac{\cos \btheta - \alpha}{1-\alpha \cos \btheta}.
\end{align*}
Because of the assumption $v^{-2}(\pi) > v^{-2}(\frac{\pi}{2}^+)$, there exist a positive number $1>\epsilon>0$ such that
\begin{align*}
v^{-2}\big(\frac{\pi}{2} + \epsilon\big) < (1-\epsilon) v^{-2}(\pi).
\end{align*} 
Choose $\alpha_0>0$ sufficiently small such that
\begin{align*}
\cos \big(\frac{\pi}{2} + \epsilon \big) < -\alpha_0,
\quad
\sqrt{\frac{1+\alpha_0}{1-\alpha_0}} \cdot (1-\epsilon) <1.
\end{align*}
Then for any $\alpha \in (0, \alpha_0]$, we have
\begin{align*}
\cos \big(\frac{\pi}{2} + \epsilon \big) < -\alpha_0 \leq -\alpha,
\quad
\sqrt{\frac{1+\alpha}{1-\alpha}} \cdot (1-\epsilon) \leq \sqrt{\frac{1+\alpha_0}{1-\alpha_0}} \cdot (1-\epsilon) <1.
\end{align*}
The first inequality implies that if $\btheta \in [0,\frac{\pi}{2}]$, then $\theta \in [0, \frac{\pi}{2} + \epsilon]$. We show that the above chosen $\alpha$ implies that 
\begin{align*}
\max_{\theta} \{ \varphi_{\alpha, 0, 0} (v^{-2}) (\theta) \} < \max_{\theta} \{ v^{-2} (\theta) \} = v^{-2}(\pi).
\end{align*}
Since for $\btheta \in [\frac{\pi}{2}, \pi]$,
\begin{align*}
\varphi_{\alpha,0,0}(v^{-2}) (\btheta) = \frac{\sqrt{1-\alpha^2}}{1 - \alpha \cos \btheta} v^{-2}(\theta) \leq \sqrt{1-\alpha^2} v^{-2}(\pi) < v^{-2}(\pi),
\end{align*}
and for $\btheta \in [0, \frac{\pi}{2}]$,
\begin{align*}
\varphi_{\alpha,0,0}(v^{-2}) (\btheta) 
&
= \frac{\sqrt{1-\alpha^2}}{1 - \alpha \cos \btheta} \cdot v^{-2}(\theta) 
\\
&
\leq \sqrt{\frac{1+\alpha}{1-\alpha}} \cdot v^{-2}\big(\frac{\pi}{2} + \epsilon\big) 
< \sqrt{\frac{1+\alpha}{1-\alpha}} \cdot (1- \epsilon) v^{-2} (\pi ) 
< v^{-2}(\pi).
\end{align*}
Then we prove the lemma.
\end{proof}

By the similar method, we can prove the following lemma.
\begin{lemma}\label{lem 6.2}
Let $v$ be a positive function on the circle satisfying $v(\theta) = v(-\theta)$. Let $m=\sup_{\theta} \{ v^{-2}(\theta) \}$.
\begin{enumerate}[label=\alph*.]
\item
If $\varlimsup_{x\rightarrow \theta} v_0^{-2} (x)  < m$ for every $\theta\in[\frac{\pi}{2},\pi]$, then there exists $\alpha_0 \in (0,1)$ such that for every Lorentz transformation $\varphi_{-\alpha, 0, 0}, \alpha\in (0,\alpha_0]$ we have that
\begin{align*}
\sup_{\theta} \{ \varphi_{-\alpha, 0, 0} (v^{-2}) (\theta) \} < \sup_{\theta} \{ v^{-2} (\theta) \} =m,
\quad
\inf_{\theta} \{ \gamma_{-\alpha, 0, 0} (v) (\theta) \} > \inf_{\theta} \{ v (\theta) \}.
\end{align*}
\item
If $\varlimsup_{x\rightarrow \theta} v^{-2} (x)  < m$ for every $\theta\in[0, \frac{\pi}{2}]$, then there exists $\alpha_0 \in (0,1)$ such that for every Lorentz transformation $\varphi_{\alpha, 0, 0}, \alpha\in (0,\alpha_0]$ we have that
\begin{align*}
\sup_{\theta} \{ \varphi_{\alpha, 0, 0} (v^{-2}) (\theta) \} < \sup_{\theta} \{ v^{-2} (\theta) \} =m,
\quad
\inf_{\theta} \{ \gamma_{\alpha, 0, 0} (v) (\theta) \} > \inf_{\theta} \{ v (\theta) \}.
\end{align*}
\end{enumerate}
\end{lemma}
\begin{proof}
\begin{enumerate}[label=\alph*., leftmargin=.2in]
\item
There exists $\epsilon>0$ such that $v_0^{-2} (\theta)  < (1-\epsilon) m, \theta\in[\frac{\pi}{2}-\epsilon,\pi]$. Choose $\alpha_0$ such that
\begin{align*}
\cos \big(\frac{\pi}{2} - \epsilon\big) > \alpha_0, 
\quad
\sqrt{\frac{1+\alpha_0}{1-\alpha_0}} (1-\epsilon) <1.
\end{align*}
Then the rest of the proof proceeds the same as the proof of lemma \ref{lem 6.1}. See figure \ref{fig 6} for $\varphi_{-\alpha_0,0,0}$.
\begin{figure}[H]
\centering
\begin{tikzpicture}[scale=1]
\draw (-5,5) node[above left] {$\theta=\pi$} -- (0,0) -- (5,5) node[above right] {$\theta=0$};
\draw[line width=0.2pt,domain=-1:1] plot (\x, {sqrt((2.85)^2+(\x)^2)});
\draw (-2.85,2.85) node[left] {\footnotesize $(1-\epsilon) m$} -- (0,2.85) (-3.5,3.5) node[left] {\footnotesize $m$} -- (3.5,3.5);
\draw (-0.5, 2.894) -- (-3.3,3.3);
\draw (0,0) -- (0,2.85);
\draw (0,0) -- (-0.5, 2.894);
\draw[->,line width=1pt,domain=0:-0.5] plot (\x, {sqrt((2.85)^2+(\x)^2)});
\draw[->,line width=1pt] (-2.85,2.85) -- (-3.3,3.3);
\draw[->,line width=1pt] (3.5,3.5) -- (2.5,2.5);
\draw[->,line width=1pt,domain=2.5:1.5] plot (\x, {sqrt((2.4)^2+(\x)^2)});
\end{tikzpicture}
\caption{Lorentzian transformation $\varphi_{-\alpha_0,0,0}$.}
\label{fig 6}
\end{figure}
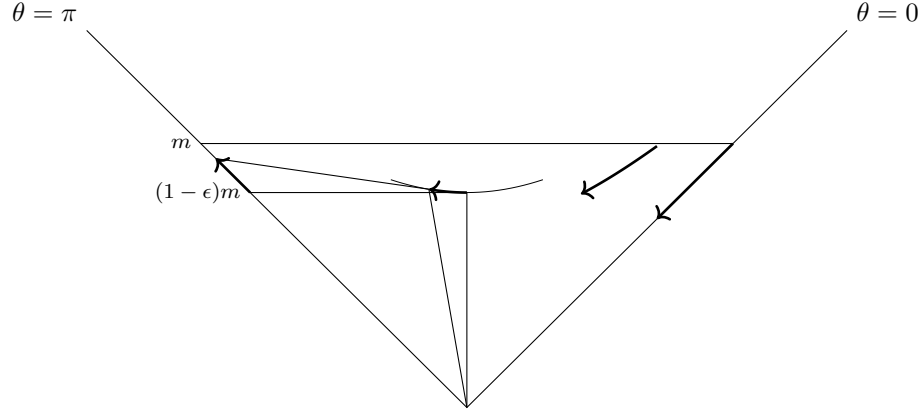
\item
There exists $\epsilon>0$ such that $v_0^{-2} (\theta)  < (1-\epsilon) m, \theta\in[0, \frac{\pi}{2}+\epsilon]$. Choose $\alpha_0$ such that
\begin{align*}
\cos \big(\frac{\pi}{2} + \epsilon\big) < - \alpha_0, 
\quad
\sqrt{\frac{1+\alpha_0}{1-\alpha_0}} (1-\epsilon) <1.
\end{align*}
Then the rest of the proof also proceeds the same as the proof of lemma \ref{lem 6.1}.
\end{enumerate}
\end{proof}

\begin{remark}\label{rem 6.3}
We have analogous lemmas of lemmas \ref{lem 6.1}, \ref{lem 6.2} for the maximum/supremum of a function and the Lorentz transformation.
\end{remark}

\subsection{Construct decreasing sequence $\{ v_n \}_{n\in \mathbb{N}}$ of $\calF$ with nondecreasing minimum}\label{subsec 6.2}
In this subsection, we use the symmetric decreasing rearrangement and Lorentz transformation to construct a decreasing sequence of $\calF$ with nondecreasing minimum starting from an arbitrary positive function $v$ on the circle. Without loss of generality, we assume $v= v^*$ since the symmetric decreasing rearrangement doesnot increase $\calF$. The idea to construct a decreasing sequence is to apply the symmetric decreasing rearrangement and Lorentz transformation alternatively.
\begin{definition}\label{def 6.4}
Let $v_0$ be a positive lower semicontinuous function on the circle satisfying $v_0(\theta) = v_0(-\theta)$ and monotonically nonincreasing on $[0,\pi]$, i.e. $v_0=v_0^*$. We construct the sequence $\{ v_n \}_{n\in \mathbb{N}}$ where $v_n = v_n^*$: suppose $v_n$ is obtained, we construct $v_{n+1}$ as follows,
\begin{enumerate}[label=\alph*.]
\item Consider the set of function $\{\varphi_{\alpha, 0, 0} (v_n^{-2})\}_{\alpha\geq 0}$. Choose $\alpha_n$ such that
\begin{align*}
\max_{\theta} \{\varphi_{\alpha_n,0,0}(v_n^{-2})(\theta)\} 
= 
\min_{\alpha\geq 0} \max_{\theta} \{\varphi_{\alpha,0,0}(v_n^{-2})(\theta)\},
\end{align*}
and for all $\alpha \in [0, \alpha_n)$,
\begin{align*}
\max_{\theta} \{\varphi_{\alpha,0,0}(v_n^{-2})(\theta)\} 
>
\min_{\alpha\geq 0} \max_{\theta} \{\varphi_{\alpha,0,0}(v_n^{-2})(\theta)\}. \footnotemark
\end{align*}
\footnotetext{This condition cannot be omitted, because of functions like $v^{-2}(\theta) = \min\{ 1, \nu_{-\alpha, 0} (\theta) \}$}.
\item Let $v_{n+1}$ be the symmetric decreasing rearrangement of $\gamma_{\alpha_n,0,0}(v_n)$.
\end{enumerate}
\end{definition}
We can easily obtain the following property of the above constructed sequence $\{v_n\}_{n\in \mathbb{N}}$.
\begin{lemma}\label{lem 6.5}
The minimum of $\{v_n(\theta)\}_{n\in \mathbb{N}}$ is monotonically nondecreasing, the maximum of $\{v_n^{-2} (\theta) \}_{n\in \mathbb{N}}$ and the functional $\calF$ of the sequence $\{ v_n \}_{n \in \mathbb{N}}$ are both monotonically nonincreasing, i.e.
\begin{align*}
\min_{\theta} \{ v_{n+1} (\theta) \} \geq \min_{\theta} \{ v_n (\theta) \},
\quad
\max_{\theta} \{ v_{n+1}^{-2} (\theta) \} \leq \max_{\theta} \{ v_n^{-2} (\theta) \},
\quad
\calF[v_{n+1}] \leq \calF[v_n],
\end{align*}
and the integral $\int_0^{2\pi} v_n^{-2}\, \ed \theta$ remains constant.
\end{lemma}
Since the sequence $\{v_n^{-2}\}_{n\in \mathbb{N}}$ is bounded and each $v_n^{-1}$ is monotonically nondecreasing on $[0,\pi]$, there exists a converging subsequence $\{v_{n_i}^{-1}\}_{i\in \mathbb{N}}$. We shall study the limit function in the next subsection. We prove a useful lemma on the step \textit{a.} in definition \ref{def 6.4}, which plays an important role in the next subsection.
\begin{lemma}\label{lem 6.6}
Let $v$ be a positive lower semicontinuous function on the circle satisfying $v(\theta) = v(-\theta)$ and monotonically nonincreasing on $[0,\pi]$, i.e. $v=v^*$. Let $\alpha_v$ chosen as in definition \ref{def 6.4}.a.,
\begin{align*}
\max_{\theta} \{\varphi_{\alpha_v,0,0}(v^{-2})(\theta)\} 
= 
\min_{\alpha\geq 0} \max_{\theta} \{\varphi_{\alpha,0,0}(v^{-2})(\theta)\}.
\end{align*}
Then either of the two assertions is ture.
\begin{enumerate}[label=\alph*.]
\item $\varphi_{\alpha_v,0,0}(v^{-2})(\frac{\pi}{2}) = \max_{\theta} \{\varphi_{\alpha_v,0,0}(v^{-2})(\theta)\}$.
\item there exist $\theta_1 \in [0, \frac{\pi}{2}) $ and $\theta_2 \in (\frac{\pi}{2}, \pi]$ such that
\begin{align*}
\varphi_{\alpha_v,0,0}(v^{-2})(\theta_1) = \varphi_{\alpha_v,0,0}(v^{-2})(\theta_2) = \max_{\theta} \{\varphi_{\alpha_v,0,0}(v^{-2})(\theta)\}.
\end{align*}
\end{enumerate}
\end{lemma}
\begin{proof}
We prove the lemma by the argument of contradictions. Let $ m = \max_{\theta} \{\varphi_{\alpha_v,0,0}(v^{-2})(\theta)\}$. Suppose that assertions \textit{a.} and \textit{b.} are both wrong. Let $v_0 = \gamma_{\alpha_v,0,0}(v)$, i.e. $v_0^{-2} = \varphi_{\alpha_v,0,0}(v^{-2})$. Then there are two cases:
\begin{enumerate}[label=\roman*.]
\item
$\varlimsup_{x\rightarrow \theta} v_0^{-2} (x) =v_0^{-2} (\theta)  < m, \forall \theta\in[\frac{\pi}{2},\pi]$,
\item
$\varlimsup_{x\rightarrow \theta} v_0^{-2} (x) = v_0^{-2}(\theta) < m, \forall \theta\in[0, \frac{\pi}{2}]$.
\end{enumerate}

In case i. by lemma \ref{lem 6.6}.\textit{a.}, we have that $\alpha_v>0$ and there exists $\alpha_v > \alpha_0>0$ sufficiently small, such that 
\begin{align*}
\max_{\theta} \{ \varphi_{-\alpha_0,0,0}(v_0^{-2})(\theta) \} 
= 
\max_{\theta} \{ \varphi_{-\alpha_0,0,0}\circ \varphi_{\alpha_v,0,0}(v^{-2})(\theta) \} 
< m,
\end{align*}
which contradicts to $m = \min_{\alpha\geq 0} \max_{\theta} \{\varphi_{\alpha,0,0}(v^{-2})(\theta)\}$.

In case ii. by lemma \ref{lem 6.6}.\textit{b.}, we have that there exists $\alpha_0>0$ sufficiently small, such that 
\begin{align*}
\max_{\theta} \{ \varphi_{\alpha_0,0,0}(v_0^{-2})(\theta) \} 
= 
\max_{\theta} \{ \varphi_{\alpha_0,0,0}\circ \varphi_{\alpha_v,0,0}(v^{-2})(\theta) \} 
< m,
\end{align*}
which contradicts to $m = \min_{\alpha\geq 0} \max_{\theta} \{\varphi_{\alpha,0,0}(v^{-2})(\theta)\}$.
\end{proof}

\begin{remark}\label{rem 6.7}
We use the Lorentz transformation to reduce the maximum of the decreasing sequence of $\calF$ in definition \ref{def 6.4}. As mentioned in remark \ref{rem 6.3}, one can also use the Lorentz transformation to reduce the maximum/supremum of a function, thus one can introduce an analogous decreasing sequence of $\calF$ but with nonincreasing maximum.
\end{remark}

\subsection{Subsequence limit of decreasing sequence $\{v_n\}_{n\in \mathbb{N}}$ of $\calF$}\label{subsec 6.3}
Let $v_0$ and the decreasing sequence $\{v_n\}_{n\in \mathbb{N}}$ be as in definition \ref{def 6.4}. Let $v_{\infty}^{-1}$ be the limit of a converging subsequence $\{v_{n_i}^{-1}\}_{i\in \mathbb{N}}$, where $v_{\infty}$ takes value in the extended positive real number set $(0, +\infty]$. We also call that $\{v_{n_i}\}_{i\in \mathbb{N}}$ converges to $v_{\infty}$. $\{v_{n_i}\}_{i\in \mathbb{N}}$ and $v_{\infty}$ satisfy the following property.
\begin{lemma}\label{lem 6.8}
Let $v_0$ and the decreasing sequence $\{v_n\}_{n\in \mathbb{N}}$ be as in definition \ref{def 6.4}. Let $v_{\infty}$ be the limit of a subsequence $\{v_{n_i}\}_{i\in \mathbb{N}}$, which always exists. Then the subsequence limit $v_{\infty}$ satisfies that
\begin{enumerate}[label=\alph*.]
\item $v_{\infty}$ is monotonically nonincreasing on $[0,\pi]$ and $v_{\infty} (\theta) = v_{\infty}(-\theta)$.

\item $\max_{\theta} \{v_{\infty}^{-2}(\theta)\} \leq \max_{\theta} \{v_n^{-2}(\theta)\}, \quad
\min_{\theta} \{v_{\infty}(\theta)\} \geq \min_{\theta} \{v_n(\theta)\}, \quad \forall n\in \mathbb{N}$.

\item $\int_0^{2\pi} v_{\infty}^{-2}\, \ed \theta = \lim_{i\rightarrow \infty} \int_0^{2\pi} v_{n_i}^{-2}\, \ed \theta = \int_0^{2\pi} v_n^{-2}\, \ed \theta, \quad n \in \mathbb{N}$.

\item $v_{\infty}$ is constant on $(\frac{\pi}{2}, \pi]$. See figure \ref{fig 7}.
\begin{figure}[H]
\centering
\begin{tikzpicture}
\draw[->] (-.1*pi,0) -- (2.1*pi,0);
\draw[->] (1*pi,0) -- (1*pi,2) node[right] {$v^*$};
\draw[red] plot [smooth] coordinates 
{(0*pi,2-1.3) ( .5*pi, 2-1.3) ( .6*pi, 2-1.3) ( .7*pi, 2-1.2) ( .8*pi, 2-0.7) ( .9*pi, 2-0.5) ( 1*pi, 2-0.4)   (1.1*pi, 2-0.5) ( 1.2*pi, 2-0.7) ( 1.3*pi, 2-1.2) ( 1.4*pi, 2-1.3) ( 1.5*pi, 2-1.3) ( 2*pi, 2-1.3) };
\draw[dashed,red] (0.6*pi,2-1.3) -- (0.6*pi,0) (1.4*pi,2-1.3) -- (1.4*pi,0);
\draw (0*pi,.1) -- (0*pi,0) node[below] {$-\pi$};
\draw (.5*pi,.1) -- (.5*pi,0) node[below] {$-\frac{\pi}{2}$};
\draw (1*pi,.1) -- (1*pi,0) node[below] {$0$};
\draw (1.5*pi,.1) -- (1.5*pi,0) node[below] {$\frac{\pi}{2}$};
\draw (2*pi,.1) -- (2*pi,0) node[below] {$\pi$};
\end{tikzpicture}
\caption{$v_{\infty}$.}
\label{fig 7}
\end{figure}
\end{enumerate}
\end{lemma}
\begin{proof}
\begin{enumerate}[label=\alph*.,leftmargin=.2in]
\item 
It follows from that $v_{n_i} = v_{n_i}^*$.

\item 
It follows from lemma \ref{lem 6.5}.

\item 
$\int_0^{2\pi} v_{\infty}^{-2}\, \ed \theta = \lim_{i\rightarrow \infty} \int_0^{2\pi} v_{n_i}^{-2} \, \ed \theta$ follows from the uniform upper boundedness of $\{v_{n_i}^{-2}\}_{i \in \mathbb{N}}$. 

\item 
Proof by contradiction: assume that $v_{\infty}$ isnot constant on $(\frac{\pi}{2}, \pi]$, then $v_{\infty}^{-2}(\pi) > v_{\infty}^{-2}(\frac{\pi}{2}^+)$. Thus by lemma \ref{lem 6.1} there exists $\alpha_{\infty}>0$ such that
\begin{align*}
\max_{\theta} \{ \varphi_{\alpha_{\infty}, 0, 0} (v_{\infty}^{-2}) (\theta) \} < \max_{\theta} \{ v_{\infty}^{-2} (\theta) \}.
\end{align*}
Since $v_{n_i}$ converges to $v_{\infty}$, we have that
\begin{align*}
\max_{\theta} \{ \varphi_{\alpha_{\infty}, 0, 0} (v_{n_i}^{-2}) (\theta) \} \longrightarrow \max_{\theta} \{ \varphi_{\alpha_{\infty}, 0, 0} (v_{\infty}^{-2}) (\theta) \}
\text{ as }
i \rightarrow \infty, 
\end{align*}
thus there exists $i\in \mathbb{N}$ such that
\begin{align*}
\max_{\theta} \{ \varphi_{\alpha_{\infty}, 0, 0} (v_{n_i}^{-2}) (\theta) \} < \max_{\theta} \{ v_{\infty}^{-2}(\theta) \}.
\end{align*}
Then by definition \ref{def 6.4},
\begin{align*}
\max_{\theta} \{ v_{n_i+1}^{-2} (\theta) \} 
= 
\min_{\alpha \geq 0} \max_{\theta} \{ \varphi_{\alpha, 0, 0} (v_{n_i}^{-2}) (\theta) \} 
\leq
\max_{\theta} \{ \varphi_{\alpha_{\infty}, 0, 0} (v_{n_i}^{-2}) (\theta) \}
< 
\max_{\theta} \{ v_{\infty}^{-2}(\theta) \},
\end{align*}
which contradicts with \textit{b}. Thus $v_{\infty}$ is constant on $(\frac{\pi}{2}, \pi]$.
\end{enumerate}
\end{proof}

\subsection{Boundedness of subsequence limit of $\{ v_n \}_{n\in \mathbb{N}}$}\label{subsec 6.4}
Suppose that $v_0$ satisfies the constraint $\int_0^{2\pi} v_0^{-2} \, \ed \theta = 2\pi$. Then by above lemma \ref{lem 6.8}.\textit{c}., $v_{\infty}$ also satisfies the constraint $\int_0^{2\pi} v_{\infty}^{-2} \, \ed \theta = 2\pi$. Since $\{\calF[v_n]\}$ is nonincreasing, we wonder whether $\calF[v_{\infty}] \leq \lim_{n\rightarrow +\infty} \calF[v_n]$. To address this question, we first prove the following lemmas which ensure that $\{ v_n \}_{n \in \mathbb{N}}$ is bounded from above.

\begin{lemma}\label{lem 6.9}
Let $v_0$ and the decreasing sequence $\{v_n\}_{n\in \mathbb{N}}$ be as in definition \ref{def 6.4}. Let $\{\alpha_n\geq 0\}_{n\in \mathbb{N}}$ also be as in definition \ref{def 6.4}. Then there exists a positive number $1>\balpha>0$ such that $0\leq \alpha_n \leq \balpha$. 
\end{lemma}
\begin{proof}
Let $m_n= \max_{\theta} \{ v_n^{-2} (\theta)\}$ and $m_{\infty}= \lim_{n\rightarrow \infty} \max_{\theta} \{ v_n^{-2} (\theta)\}$. Then for any $\epsilon>0$, there exists $N_{\epsilon}$ sufficiently large such that
\begin{align*}
m_{\infty} \leq \max_{\theta} \{ v_n^{-2} (\theta)\} < m_{\infty} + \epsilon, \quad n\geq N_{\epsilon}.
\end{align*}
By lemma \ref{lem 6.6}, there exists $\btheta_{n} \in [\frac{\pi}{2}, \pi]$ that
\begin{align*}
\varphi_{\alpha_{n},0,0} (v_{n}^{-2}) (\btheta_n)
=
\max_{\theta} \{ \varphi_{\alpha_{n},0,0} (v_{n}^{-2}) (\btheta)  \} .
\end{align*}
Let $\theta_{n}\in [\btheta_{n_k},\pi]$ be defined by equations \eqref{eqn 3.1} that
\begin{align}\label{eqn 5.1}
\cos \theta_{n} 
= \frac{\cos \btheta_{n} - \alpha_{n}}{1-\alpha_{n} \cos \btheta_{n}},
\end{align}
We prove the lemma by the argument of contradictions. Suppose on the contrary, there exists a subsequence $\{\alpha_{n_k}\}_{k\in \mathbb{N}}$ converges to $1$ as $k\rightarrow +\infty$,
then we have that
\begin{align*}
&
\cos \theta_{n_k}
= -\alpha_{n_k} + (1-\alpha_{n_k}^2) \frac{\cos \btheta_{n_k}}{1-\alpha_{n_k} \cos \btheta_{n_k}}
\rightarrow
-1
\text{ as }
k \rightarrow +\infty,
\\
&
\theta_{n_k} \rightarrow \pi
\text{ as }
k \rightarrow +\infty.
\end{align*}
See figure \ref{fig 8}. We shall prove that $\int_0^{\pi} v_{n_k}^{-2}\, \ed \theta$ converges to $0$, which contradicts with $\int_0^{\pi} v_{n_k}^{-2}\, \ed \theta = \pi$. For the sake of brevity, use $\bar{v}_n$ to denote $\varphi_{\alpha_n,0,0} (v_{n})$ in the following.
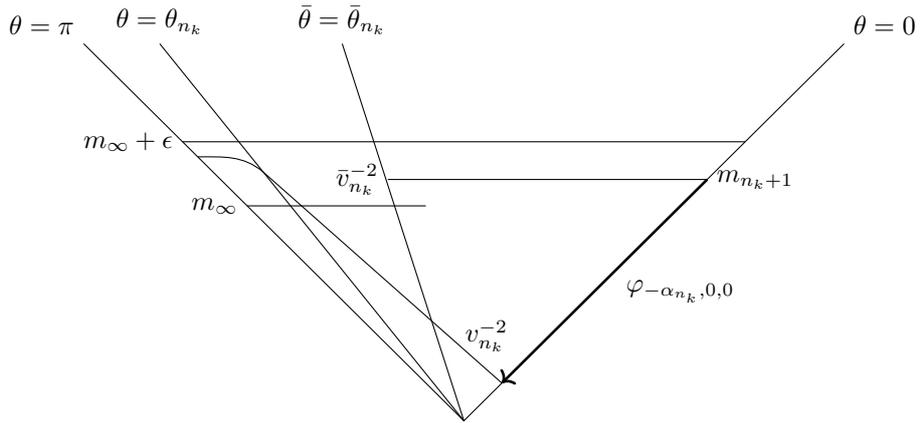
\begin{figure}[H]
\centering
\begin{tikzpicture}[scale=1]
\draw (-5,5) node[above left] {$\theta=\pi$} -- (0,0) -- (5,5) node[above right] {$\theta=0$};
\draw (0,0) -- (-1.6,5) node[above] {$\btheta = \btheta_{n_k}$};
\draw (0,0) -- (-4,5) node[above] {$\theta=\theta_{n_k}$};
\draw (-2.85,2.85) node[left] {$m_{\infty}$}  -- (-0.5,2.85) 
(-3.7,3.7) node[left] {$m_{\infty} + \epsilon$} -- (3.7,3.7)
(3.2,3.2) node[right] {$m_{n_k+1}$} -- (-1,3.2) node[left] {$\bar{v}_{n_k}^{-2}$};
\draw (-2.64,3.3) 
to [out=180-40,in=0] (-3.5,3.5);
\draw (-2.64,3.3) -- (0.5,0.5);
\node[above] at (0.3,0.8) {$v_{n_k}^{-2}$};
\draw[->, line width=1pt] (3.2,3.2) -- (2,2) node[below right] {$\varphi_{-\alpha_{n_k},0,0}$} -- (0.5,0.5);
\end{tikzpicture}
\caption{$v_{n_k}^{-2}$, $\alpha_{n_k} \rightarrow 1$, $\theta_{n_k} \rightarrow \pi$.}
\label{fig 8}
\end{figure}

On $[\theta_{n_k}, \pi]$, $v_{n_k}^{-2} < m_{\infty}+\epsilon$. On $[0, \theta_{n_k}]$,
\begin{align*}
v_{n_k}^{-2} (\theta)
&
=
\varphi_{-\alpha_{n_k},0,0} (\bar{v}_n^{-2}) (\theta)
=
\frac{\sqrt{1-\alpha_{n_k}^2}}{1+ \alpha_{n_k} \cos \theta} \cdot \bar{v}_{n_k}^{-2} (\btheta)
\leq
\frac{\sqrt{1-\alpha_{n_k}^2}}{1+ \alpha_{n_k} \cos \theta} \cdot m_{n_k+1}.
\end{align*}
Note that at $\theta_{n_k}$ we have
\begin{align*}
v_{n_k}^{-2} (\theta_{n_k})
&
=
\frac{\sqrt{1-\alpha_{n_k}^2}}{1+ \alpha_{n_k} \cos \theta_{n_k}} \cdot \bar{v}_{n_k}^{-2} (\btheta_{n_k})
=
\frac{\sqrt{1-\alpha_{n_k}^2}}{1+ \alpha_{n_k} \cos \theta_{n_k}} \cdot m_{n_k+1} < m_{\infty}+\epsilon,
\quad
\text{ for }
n_k \geq N_{\epsilon}.
\end{align*}
then on $[0,\theta_{n_k}]$
\begin{align*}
v_{n_k}^{-2} (\theta) 
&
\leq 
(m_{\infty}+\epsilon) \cdot \frac{1+ \alpha_{n_k} \cos \theta_{n_k}}{1+ \alpha_{n_k} \cos \theta}
\\
&
\overset{\eqref{eqn 5.1}}{=}
(m_{\infty}+\epsilon) \cdot \frac{\sqrt{1-\alpha_{n_k}^2}}{1-\alpha_{n_k} \cos \btheta_{n_k}} \cdot \frac{\sqrt{1-\alpha_{n_k}^2}}{1+ \alpha_{n_k} \cos \theta}
\\
&
\leq
(m_{\infty}+\epsilon) \sqrt{1-\alpha_{n_k}^2} \cdot \nu_{\alpha_{n_k},0} (\theta),
\end{align*}
where the last inequality follows from $\alpha_{n_k} \geq 0$, $\btheta_{n_k} \in [\frac{\pi}{2}, \pi]$. Hence we can estimate $\int_0^{\pi} v_{n_k}^{-2}\, \ed \theta$ now,
\begin{align*}
\int_0^{\pi} v_{n_k}^{-2}\, \ed \theta 
&=
\int_0^{\theta_{n_k}} v_{n_k}^{-2}\, \ed \theta 
+ 
\int_{\theta_{n_k}}^{\pi} v_{n_k}^{-2}\, \ed \theta
\\
&\leq
\int_0^{\theta_{n_k}} (m_{\infty}+\epsilon) \sqrt{1-\alpha_{n_k}^2} \cdot \nu_{\alpha_{n_k},0} (\theta) \, \ed \theta 
+ 
\int_{\theta_{n_k}}^{\pi} (m_{\infty}+\epsilon) \, \ed \theta
\\
&\leq
(m_{\infty}+\epsilon) \sqrt{1-\alpha_{n_k}^2} \pi + (m_{\infty}+\epsilon) (\pi - \theta_{n_k})
\\
&
\rightarrow
0 
\text{ as }
k \rightarrow + \infty.
\end{align*}
Thus we arrive at the desired contradiction.
\end{proof}

In fact, we have a much stronger result, which implies lemma \ref{lem 6.9} as a corollary.
\begin{lemma}\label{lem 6.10}
Let $v_0$ and the decreasing sequence $\{v_n\}_{n\in \mathbb{N}}$ be as in definition \ref{def 6.4}. Let $\{\alpha_n\geq 0\}_{n\in \mathbb{N}}$ also be as in definition \ref{def 6.4}. Then 
\begin{align*}
\prod_{n\in \mathbb{N}} \sqrt{1-\alpha_n^2} \geq \frac{m_{\infty}}{m_0}.
\end{align*}
\end{lemma}
\begin{proof}
Adopt the notations in lemma \ref{lem 6.9}. We have that
\begin{align*}
\cos \theta_n
= 
\frac{\cos \btheta_{n} - \alpha_{n}}{1-\alpha_{n} \cos \btheta_{n}} \overset{\btheta_n \in [\frac{\pi}{2}, \pi]}{\leq} -\alpha_n.
\end{align*} 
By the formula
\begin{align*}
\varphi_{\alpha_{n},0,0} (v_{n}^{-2}) (\btheta_n) 
= 
\frac{1+\alpha_n \cos\theta_n}{\sqrt{1-\alpha_n^2}} v_{n}^{-2} (\theta_n)
=
\max_{\theta} \{ \varphi_{\alpha_{n},0,0} (v_{n}^{-2}) (\btheta)  \}
=
\max_{\theta} \{ v_{n+1}^{-2} (\btheta)  \},
\end{align*}
we obtain that
\begin{align*}
\sqrt{1-\alpha_n^2} \cdot m_n \geq m_{n+1}.
\end{align*}
Thus the lemma follows easily by taking the products of the above inequality.
\end{proof}

\begin{lemma}\label{lem 6.11}
Let $v_0$ and the decreasing sequence $\{v_n\}_{n\in \mathbb{N}}$ be as in definition \ref{def 6.4}. Then $\{ v_n\}_{i \in \mathbb{N}}$ is bounded by a finite positive number from above.
\end{lemma}
\begin{proof}
It is sufficient to prove that $\{v_n^{-2}\}_{n\in\mathbb{N}}$ is bounded by a positive number from below. Adopt the notations in lemma \ref{lem 6.9}, \ref{lem 6.10}.  We estimate the lower bound of $\varphi_{\alpha_{n},0,0} (v_{n}^{-2}) (\btheta)$. Firstly, we define $\btheta'_n\in [0,\frac{\pi}{2}], \theta'_n$ similarly as $\btheta_n, \theta_n$ in lemma \ref{lem 6.9} by
\begin{align*}
\varphi_{\alpha_{n},0,0} (v_{n}^{-2}) (\btheta'_n)
=
\max_{\theta} \{ \varphi_{\alpha_{n},0,0} (v_{n}^{-2}) (\btheta)  \},
\quad
\cos \theta'_{n} 
= \frac{\cos \btheta'_{n} - \alpha_{n}}{1-\alpha_{n} \cos \btheta'_{n}}.
\end{align*}
See figure \ref{fig 9}. The existence of such $\btheta'_n$ is confirmed by lemma \ref{lem 6.6}.
\begin{figure}[H]
\centering
\begin{tikzpicture}[scale=1]
\draw (-5,5) node[above left] {$\theta=\pi$} -- (-0,0) -- (5,5) node[above right] {$\theta=0$};
\draw (0,0) -- (1.6,5) node[above] {\footnotesize $\btheta'_{n}$};
\draw (0,0) -- (0,5) node[above] {$\btheta=\frac{\pi}{2}$};
\draw (0,0) -- (-1.6,5) node[above] {\footnotesize $\theta'_{n}$};
\draw (0,0) -- (-3,5) node[above] {\footnotesize $\cos \theta =-\alpha_n $};
\draw[->,line width =.5pt,domain=-1.6*0.74:1.6*0.74] plot (\x, {sqrt(3.7^2+(\x)^2)});
\draw[<-,line width =.5pt,domain=0:-3*0.87] plot (\x, {sqrt(3.6^2+(\x)^2)});
\draw[<-,domain=0:-3*0.87*0.3] plot (\x, {sqrt((3.5*0.3)^2+(\x)^2)});
\draw (-3*0.87*0.3,5*0.87*0.3) -- (5*0.87*0.3, 5*0.87*0.3) node[below right] {$\inf_{\theta} \{ v_n^{-2} (\theta) \}$};
\draw (0,3.5*0.3) -- (3.5*0.7,3.5*0.7);
\draw[->,line width=1pt] (5*0.87*0.3, 5*0.87*0.3) -- (3.5*0.5,3.5*0.5) node[right] {$\varphi_{\alpha_n,0,0}$} -- (3.5*0.7,3.5*0.7);
\draw (5*0.74-0.5, 5*0.74-0.5) node[right] {$m_{\infty}$} -- (1*0.74,5*0.74-0.5);
\draw (1.6*0.74,5*0.74) node[right] {$m_{n+1}$} -- (1.6*0.74-3.5*0.7*1.26, 5*0.74 +3.5*0.3*1.26-3.5*0.7*1.26);
\draw (-1.6*0.74,5*0.74) -- (-5*0.74,5*0.74) node[left] {$v_n^{-2}(\theta'_n)$};
\draw[->,line width=1pt] (-5*0.74,5*0.74) -- (-4*0.74,4*0.74) node[below left] {$\varphi_{\alpha_n,0,0}$} -- (1.6*0.74- 3.5*0.7*1.26, 5*0.74 +3.5*0.3*1.26-3.5*0.7*1.26);
\end{tikzpicture}
\caption{Lower bound of $\varphi_{\alpha_n,0,0}(v_n^{-2})$.}
\label{fig 9}
\end{figure}

For $\btheta\in [0, \frac{\pi}{2}]$,
\begin{align*}
\varphi_{\alpha_{n},0,0} (v_{n}^{-2}) (\btheta) 
= 
\frac{\sqrt{1-\alpha^2}}{1 - \alpha \cos \btheta}  \cdot v_{n}^{-2} (\theta)
&\geq
\sqrt{1-\alpha_n^2} \cdot v_n^{-2} (\theta), \quad \btheta \in [0, \frac{\pi}{2}],
\\
&\geq
\sqrt{1-\alpha_n^2} \cdot \inf_{\theta} \{ v_n^{-2} (\theta) \}.
\end{align*}
For $\btheta \in [\btheta'_n, \pi]$, which corresponds to $\theta\in [\theta'_n, \pi]$,
\begin{align*}
\varphi_{\alpha_{n},0,0} (v_{n}^{-2}) (\btheta) 
&
=
\frac{\sqrt{1-\alpha_n^2}}{1 - \alpha_n \cos \btheta}  \cdot v_{n}^{-2} (\theta)
\geq
\frac{\sqrt{1-\alpha_n^2}}{1 - \alpha_n \cos \btheta}  \cdot v_{n}^{-2} (\theta'_n)
\\
&
=
\frac{\sqrt{1-\alpha_n^2}}{1 - \alpha_n \cos \btheta} \cdot \frac{1 - \alpha_n \cos \btheta'_n}{\sqrt{1-\alpha_n^2}} \cdot
\varphi_{\alpha_{n},0,0} (v_{n}^{-2}) (\btheta'_n) 
\\
&
\geq
\frac{1 - \alpha_n \cos \btheta'_n}{1 - \alpha_n \cos \btheta} \cdot m_{\infty}
\quad\quad
\big( 0 \leq \btheta'_n\leq \btheta \leq \pi \big)
\\
&
\geq
\frac{1-\alpha_n}{1 + \alpha_n} \cdot m_{\infty}
\quad\quad
(0\leq \alpha_n \leq \balpha<1 \text{ by lemma \ref{lem 6.9}})
\\
&
\geq
\frac{1-\balpha}{1+\balpha} m_{\infty}.
\end{align*}
Thus
\begin{align*}
\inf_{\theta} \{ v_{n+1}^{-2}(\theta) \}
&
=
\inf_{\btheta} \{ \varphi_{\alpha_{n},0,0} (v_{n}^{-2}) (\btheta) \}
\geq
\min \Big\{ \sqrt{1-\alpha_n^2} \cdot \inf_{\theta} \{ v_n^{-2}(\theta) \}, \frac{1-\balpha}{1+\balpha} m_{\infty}  \Big\}
\\
&
\geq
\prod_{k\in \mathbb{N}} \sqrt{1-\alpha_k^2} \cdot \frac{1-\balpha}{1+\balpha} m_{\infty}
\geq
\frac{1-\balpha}{1+\balpha} \cdot \frac{1}{m_0} >0.
\end{align*}
Then the lemma follows.
\end{proof}

\subsection{Functional $\calF$ of subsequence limit of $\{ v_n \}_{n\in \mathbb{N}}$}\label{subsec 6.5}
We have shown that the decreasing sequence $\{ v_n \}_{n\in \mathbb{N}}$ constructed in definition \ref{def 6.4} is bounded below by a positive number (lemma \ref{lem 6.8}) and bounded above by a finite positive number (lemma \ref{lem 6.11}). Then we can derive the consequence on the functional $\calF$ of the sequence $\{v_n \}_{n\in\mathbb{N}}$.
\begin{proposition}\label{pro 6.12}
Let $v_0$ and the decreasing sequence $\{v_n\}_{n\in \mathbb{N}}$ be as in definition \ref{def 6.4}. Let $v_{\infty}$ be the limit of a subsequence $\{v_{n_k} \}_{k \in \mathbb{N}}$. Suppose that $v_0 \in \mathrm{H}^1(\mathbb{S}^1)$, then we have the following assertions:
\begin{enumerate}[label=\alph*.]
\item $\{v_n \}_{n\in \mathbb{N}}$ is a bounded sequence in $\mathrm{H}^1 (\mathbb{S}^1)$,
\item $\{ v_{n_k}\}_{k\in \mathbb{N}}$ uniformly converges to $v_{\infty}$ and $v_{\infty}$ is continuous,
\item $v_{\infty}\in \mathrm{H}^1(\mathbb{S}^2)$ and $\calF[v_{\infty}] \leq \lim_{n\rightarrow + \infty}  \calF[v_n]$.
\end{enumerate}
\end{proposition}
\begin{proof} With the preparations in previous subsections, it is easy to prove the proposition now.
\begin{enumerate}[label=\alph*.]
\item 
It follows from $\{ \calF[v_n] \}_{n\in \mathbb{N}}$ is nonincreasing and $v_n$ is bounded above by a finite positive number.
\item 
\textit{a.} implies that $\{ v_n \}_{n\in \mathbb{N}}$ is equicontinuous, thus the convergence of $\{ v_{n_k}\}_{k\in \mathbb{N}}$ is uniform and the limit function $v_{\infty}$ is continuous.
\item
$\{v_{n_k} \}$ converges to $v_{\infty}$ weakly in $\mathrm{H}^1(\mathbb{S}^1)$, thus $v_{\infty} \in \mathrm{H}^1(\mathbb{S}^1)$. $\calF[v_{\infty}] \leq \lim_{n\rightarrow + \infty}  \calF[v_n]$ follows from
\begin{align*}
\int_0^{2\pi} (v'_{\infty})^2 \ed \theta \leq \lim_{k\rightarrow + \infty}  \int_0^{2\pi} (v'_{n_k})^2 \ed \theta,
\quad
\int_0^{2\pi} v_{\infty}^2 \ed \theta = \lim_{k\rightarrow + \infty}  \int_0^{2\pi} v_{n_k}^2 \ed \theta,
\quad
\calF[v_{n+1}] \leq \calF[v_n].
\end{align*}
\end{enumerate}
\end{proof}

Then we have the following proposition on the global minimum/infimum of the functional $\calF$.
\begin{proposition}\label{pro 6.13}
Let $S$ be the set of functions on the circle satisfying the following conditions:
\begin{enumerate}[label=\alph*.]
\item $v>0$ is symmetric nonincreasing on $[-\pi,\pi]$,
\item $v \in \mathrm{H}^1(\mathbb{S}^1)$,
\item $v$ is constant on $[\frac{\pi}{2}, \pi]$.
\end{enumerate}
Let $S_c$ be the subset of $S$ with the additional constraint condition:
\begin{enumerate}[label=\alph*.]\addtocounter{enumi}{3}
\item $\int_0^{2\pi} v^{-2}\, \ed \theta = 2\pi$.
\end{enumerate}
Then we have that
\begin{align*}
\inf \Big\{ \calF[v] : v\in\mathrm{H}^1(\mathbb{S}^1), v>0, \int_0^{2\pi} v^{-2}\, \ed \theta = 2\pi \Big\} = \inf_{v\in S_c} \{ \calF[v] \},
\end{align*}
and the infimum of $\calF$ on the left can be achieved if and only if that it can be achieved in $S_c$.
\end{proposition}
\begin{proof}
The proposition follows from the decreasing sequence of $\calF$ constructed in definition \ref{def 6.4}.
\end{proof}

\begin{remark}\label{rem 6.14}
There is an analogy of proposition \ref{pro 6.13} with condition \textit{c.} replaced by
\begin{enumerate}[label=\alph*\ensuremath{'}.]\addtocounter{enumi}{2}
\item $v$ is constant on $[0,\frac{\pi}{2}]$.
\end{enumerate}
It can be derived in the same way by replacing the decreasing sequence $\{ v_n \}_{n\in \mathbb{N}}$ in definition \ref{def 6.4} by the decreasing sequence of $\calF$ with nonincreasing maximum as mentioned in remark \ref{rem 6.7}.
\end{remark}

\section{Critical points are global minimisers: proof of the inequality}\label{sec 7}

In this section, we prove claim \ref{cla 4.3} that the critical points $v^{-2} = \nu_{\alpha,\theta_0}$ of the functional $\calF[v] = \int_0^{2\pi} [ 4(v')^2 - v^2 ] \ed \theta$ are actually global minimisers under the constraint $\int_0^{2\pi} v^{-2} \ed \theta =2\pi$ using proposition \ref{pro 6.13}.

Recall that proposition \ref{pro 6.13} says that in order to find the infimum of the functional $\calF[v]$ for positive functions in $\mathrm{H}^1(\mathbb{S}^1)$ under the constraint $\int_0^{2\pi} v^{-2} \ed \theta =2\pi$, it is sufficient to find the infimum of $\calF$ in a much restrictive class of functions $S_c$ where $v$ is constant on $[\frac{\pi}{2}, \pi]$. We shall show that a minimising sequence of $\calF$ in $S_c$ converges to a limit function $v_{\infty} \in S_c$, which implies that the infimum of $\calF$ can be achieved and the sequence converges to the constant function.

\begin{proof}[Proof of claim \ref{cla 4.3}]
Let $\{v_n\}$ be a minimising sequence of the functional $\calF$ in $S_c$ defined in proposition \ref{pro 6.13}. Introduce the notations
\begin{align*}
m_n = \min_{\theta \in [0,\pi]} \{ v_n(\theta) \},
\quad
M_n = \max_{\theta \in [0,\pi]} \{ v_n(\theta) \}.
\end{align*}
Since $v_n \in S_c$, we have that $v_n([\frac{\pi}{2}, \pi]) = \{m_n\}$, $v_n(0) = M_n$ and
\begin{align*}
2\pi=\int_0^{2\pi} v_n^{-2} \ed \theta 
\left\{
\begin{aligned}
&
\geq 2 \int_{\frac{\pi}{2}}^{\pi} v_n^{-2} \ed \theta = \frac{\pi}{m_n^2}
\quad
\Rightarrow
\quad
m_n \geq \frac{\sqrt{2}}{2},
\\
&
\leq 
2 \int_{0}^{\pi} m_n^{-2} \ed \theta = \frac{2\pi}{m_n^2}
\quad
\Rightarrow
\quad
m_n \leq 1.
\end{aligned}
\right.
\end{align*}

Now if $\{v_n\}$ is also uniformly bounded from above, i.e. $\sup \{ M_n \} < +\infty$, then $\{ v_n \}$ is bounded in $\mathrm{H}^1(\mathbb{S}^1)$, since
\begin{align*}
\Vert v_n \Vert_{\mathrm{H}^1(\mathbb{S}^1)} = \frac{1}{2} (\calF[v_n] + \Vert v_n \Vert_{L^2(\mathbb{S}^1)}^2)^{\frac{1}{2}}  + \Vert v_n \Vert_{L^2(\mathbb{S}^1)} 
\end{align*}
is uniformly bounded. Thus by passing to a subsequence, we can always assume that $\{v_n\}$ converges weakly to a function $v_{\infty} \in \mathrm{H}^1(\mathbb{S}^1)$, which is also a uniform convergence. Then by the boundedness of $\{ v_n \}$ and $\{ v_n^{-1} \}$, we have $v_{\infty} \in S_c$. $v_{\infty}$ achieves the global infimum of $\calF$ under the constraint as $\{v_n\}$ is a minimising sequence and
\begin{align*}
\calF[v_{\infty}] \leq \lim_{n\rightarrow +\infty} \calF[v_n].
\end{align*}
By proposition \ref{pro 4.1}, $v_{\infty}$ must be the constant function $1$. Then the claim follows.

Therefore in order to prove the claim, it is sufficient to show that $\{v_n\}$ is uniformly bounded from above. We prove this by the method of contradiction in the following. Assume the opposite that $\sup \{ M_n \} = +\infty$, then we can always assume that $M_n \rightarrow +\infty$ as $n\rightarrow +\infty$ by passing to a subsequence. We show that $\calF[v_n]$ approaches $+\infty$ as $n \rightarrow +\infty$. Note that
\begin{align*}
(M_n - m_n)^2 
&=
\Big[ \int_0^{\frac{\pi}{2}} v_n' \ed \theta \Big]^2
\leq
\frac{\pi}{2} \cdot \int_0^{\frac{\pi}{2}} (v_n')^2 \ed \theta,
\end{align*}
thus
\begin{align*}
\calF[v_n]
&= 2 \int_{0}^{\pi} [ 4(v_n')^2 -v_n^2 ] \ed \theta
\geq
\frac{16}{\pi}(M_n-m_n)^2 - \pi M_n^2 -\pi m_n^2
\rightarrow +\infty,
\end{align*}
as $M_n \rightarrow + \infty$, $\frac{\sqrt{2}}{2} \leq m_n \leq 1$. It contradicts $v_n$ being a minimising sequence. Thus we show that $\{v_n\}$ is uniformly bounded from above and the proof is complete.
\end{proof}
We see that proposition \ref{pro 6.13} plays a key role in the above proof, which restricts the class of functions in the variational problem. Now we gives the complete statement of the inequality in the following.
\begin{theorem}\label{thm 7.1}
Let $v$ be a function in $\mathrm{H}^1(\mathbb{S}^1)$, then $v$ satisfies the following inequality
\begin{align*}
\int_{\mathbb{S}^1} [4(v')^2 - v^2]  \ed \theta \geq  - \frac{4\pi^2}{\int_{\mathbb{S}^1} v^{-2} \ed \theta},
\end{align*}
with the equality being achieved at $v^{-2} (\theta) = \frac{k\sqrt{1-\alpha^2}}{1+ \alpha \cos(\theta - \theta_o)}$ where $k>0$, $\alpha \in (-1,1)$, $\theta_0 \in \mathbb{R}$.

\end{theorem}
\begin{proof}
First note that it is sufficient to prove the inequality for non-negative functions, since one can replace $v$ by $|v|$. For a non-negative function $v$, define $v_n = \max\{ v, 1/n\}$. Then the inequality of $v$ follows from taking the limit of the inequality of $v_n$ as $1/n\rightarrow 0^+$.
\end{proof}
\begin{corollary}\label{cor 7.2}
In fact if $v\in \mathrm{H}^1(\mathbb{S}^1)$ vanishes at some point, then $\int_{\mathbb{S}^1} v^{-2} \ed \theta=+\infty$,\footnote{\label{footnote 4}The proof of this fact we know is however rather cumbersome compared to the elegance of the fact. We need to study the minimum of the Dirichlet energy $E(v) = \int_0^{\pi} (v')^2 \ed \theta$ in the class of functions $B_{m,c}^M$ defined by
\begin{enumerate}[label=\alph*.]
\item $v\in \mathrm{H}^1([0,\pi])$,
\item $\min_{\theta\in [0,\pi]}\{ v(\theta) \} = m$, $\max_{\theta\in [0,\pi]}\{ v(\theta) \} =M$ and $\int_0^{\pi} v^{-2} \ed \theta = c \in (\frac{\pi}{M^2}, \frac{\pi}{m^2})$.
\end{enumerate}
Then study the limit of the minimum of the Dirichlet energy $\min_{v\in B_{m,c}^M}\{E(v)\}$ as $m\rightarrow 0^+$ while $M$, $c$ being fixed, which gives that
\begin{align*}
\min_{v\in B_{m,c}^M } \{ E(v) \} \rightarrow + \infty,
\text{ as }
m \rightarrow 0^+,
\end{align*}
where $c$, $M$, $M^{-1}$ could vary but are uniformly bounded when taking the limit.}  thus 
\begin{align*}
4\int_{\mathbb{S}^1} (v')^2 \ed \theta \geq \int_{\mathbb{S}^1} v^2  \ed \theta \text{ if } v\in \mathrm{H}^1(\mathbb{S}^1) \text{ vanishes somewhere.}
\end{align*}
Note this is similar to the Poincar\'e inequality on the circle $\int_{\mathbb{S}^1} (v')^2  \ed \theta \geq  \int_{\mathbb{S}^1} (v-\overline{v})^2  \ed \theta$ where $\overline{v}= \frac{\int_{\mathbb{S}^1} v \ed \theta}{2\pi}$.
\end{corollary}
We leave the proof of corollary \ref{cor 7.2} in appendix \ref{appen A}.

\section{Some variants of the inequality}\label{sec 8}

In this section, we derive some variants of the inequality $\int_{\mathbb{S}^1} [4(v')^2 - v^2]  \ed \theta \geq  - \frac{4\pi^2}{\int_{\mathbb{S}^1} v^{-2} \ed \theta}$.

\subsection{Inequality for functions on $[0,l]$}
It is a natural question to ask whether a similar inequality holds for a function $v$ on $[0,2\pi]$ with $v(0) \neq v(2\pi)$. We have the following corollary of theorem \ref{thm 7.1}.
\begin{corollary}\label{cor 8.1}
Let $v$ be a function in $\mathrm{H}^1([0,l])$, then $v$ satisfies the following inequality
\begin{align*}
\int_0^{l} [\frac{4l^2}{\pi^2}(v')^2 - v^2]  \ed \theta \geq  - \frac{l^2}{\int_0^{l} v^{-2} \ed \theta},
\end{align*}
with the equality achieved at $v^{-2}(\theta) = \frac{k\sqrt{1-\alpha^2}}{1+ \alpha \cos (\frac{\pi}{l}\theta)}$, $k>0$, $\alpha \in (-1,1)$.
\end{corollary}
\begin{proof}
Note that we can extend a continuous function $v$ on $[0,\pi]$ to a even periodic continuous function of period $2\pi$, thus by theorem \ref{thm 7.1} we have that
\begin{align*}
\int_0^{\pi} [4(v')^2 - v^2]  \ed \theta \geq  - \frac{\pi^2}{\int_0^{\pi} v^{-2} \ed \theta},
\end{align*}
with the equality achieved at $v^{-2}(\theta) = \frac{k\sqrt{1-\alpha^2}}{1+ \alpha \cos \theta}$,
$k>0$, $\alpha \in (-1,1)$. By a simple rescaling $\theta \rightarrow \frac{\pi}{l} \theta$, the corollary follows.
\end{proof}
\begin{corollary}\label{cor 8.2}
Following from corollary \ref{cor 7.2}, we have that
\begin{align*}
\frac{4l^2}{\pi^2} \int_0^{l} (v')^2  \ed \theta \geq  \int_0^{l} v^2  \ed \theta, 
\end{align*}
if $v\in \mathrm{H}^1([0, l])$ vanishes somewhere.
\end{corollary}

\subsection{Transformation by stereographic projection}

We transform the inequality on the circle to an inequality on the line by the stereographic projection. Recall the stereographic projection from the circle to the real line:
\begin{align*}
&
\psi:
\quad
\mathbb{S}^1 \rightarrow \mathbb{R},
\quad
\theta \mapsto x = \frac{\sin \theta}{1-\cos \theta} = \cot \frac{\theta}{2}.
\end{align*}
We have that
\begin{align*}
\ed x =  \frac{\ed \theta}{\cos \theta -1}, 
\quad
\ed \theta = \frac{2 \ed x}{1+x^2}.
\end{align*}
Thus
\begin{align*}
&
\Big| \frac{\ed v}{\ed \theta} \Big|^2 \ed \theta = \frac{1+x^2}{2}\Big| \frac{\ed v}{\ed x} \Big|^2 \ed x,
\quad
v^2 \ed \theta = \frac{2v^2}{1+x^2} \ed x,
\\
&
v^{-2} \ed \theta 
=
\Big( \sqrt{\frac{1+x^2}{2}} v \Big)^{-2} \ed x
=
\frac{2}{1+x^2} v^{-2} \ed x,
\end{align*}
Thus we obtain the following corollary.
\begin{corollary}\label{cor 8.3}
Let $v$ be a function in $\mathrm{H}^1_{loc}(\mathbb{R})$.
\begin{enumerate}[label=\alph*.]
\item If the limit $v(\infty) = \lim_{x\rightarrow +\infty} v(x) = \lim_{x\rightarrow -\infty} v(x)$ exists, then
\begin{align*}
\int_{-\infty}^{+\infty}  \Big[ (1+x^2) \Big|\frac{\ed v}{\ed x} \Big|^2 -  \frac{v^2}{1+x^2} \Big]\ed x
\geq
-\frac{\pi^2}{\int_{-\infty}^{+\infty}  \frac{v^{-2}}{1+x^2} \ed x},
\end{align*}
with the equality being achieved at $v^{-2} =  \frac{k\sqrt{1-\alpha^2}}{1+ \alpha \big[\frac{1-x_0^2}{1+x_0^2}\cdot \frac{1-x^2}{1+x^2} + \frac{2x_0}{1+x_0^2}\cdot \frac{2x}{1+x^2} \big]}$, where $k>0$, $\alpha \in (-1,1)$, $x_0 \in \mathbb{R}$.

\item If limits $v(+\infty) = \lim_{x\rightarrow +\infty} v(x)$, $ v(-\infty) = \lim_{x\rightarrow -\infty} v(x)$ both exist, then
\begin{align*}
\int_{-\infty}^{+\infty}  \Big[ 4(1+x^2) \Big|\frac{\ed v}{\ed x} \Big|^2 -  \frac{v^2}{1+x^2} \Big]\ed x
\geq
-\frac{\pi^2}{\int_{-\infty}^{+\infty}  \frac{v^{-2}}{1+x^2} \ed x},
\end{align*}
with the equality being achieved at $v^{-2} =  \frac{k\sqrt{1-\alpha^2}}{1+ \frac{\alpha x}{\sqrt{1+x^2}}}$, where $k>0$, $\alpha \in (-1,1)$.
\end{enumerate}
\end{corollary}
\begin{proof}
\textit{a.} is equivalent to $\int_{\mathbb{S}^1} [4(v')^2 - v^2]  \ed \theta \geq  - \frac{4\pi^2}{\int_{\mathbb{S}^1} v^{-2} \ed \theta}$, and \textit{b.} is equivalent to $\int_0^{2\pi} [16(v')^2 - v^2]  \ed \theta \geq  - \frac{4 \pi^2}{\int_0^{2 \pi} v^{-2} \ed \theta}$.
\end{proof}
\begin{corollary}\label{cor 8.4}
Following corollaries \ref{cor 7.2}, \ref{cor 8.2}, we have that
\begin{enumerate}
\item If $v$ vanishes somewhere on $(-\infty,+\infty) \cup \{\infty\}$, then
\begin{align*}
\int_{-\infty}^{+\infty}  (1+x^2) \Big|\frac{\ed v}{\ed x} \Big|^2 \ed x
\geq
\int_{-\infty}^{+\infty}  \frac{v^2}{1+x^2} \ed x.
\end{align*}

\item If $v$ vanishes somewhere on $[-\infty,+\infty]$, then
\begin{align*}
\int_{-\infty}^{+\infty}  4(1+x^2) \Big|\frac{\ed v}{\ed x} \Big|^2 \ed x
\geq
\int_{-\infty}^{+\infty}  \frac{v^2}{1+x^2} \ed x.
\end{align*}
\end{enumerate}
\end{corollary}

Define $u = \sqrt{\frac{1+x^2}{2}} v$ and rewrite the inequality in terms of $u$. We have the following corollary.
\begin{corollary}\label{cor 8.5}
Let $u \in \mathrm{H}^1_{loc}(\mathbb{R})$.
\begin{enumerate}[label=\alph*.]
\item If the limit $\lim_{x\rightarrow +\infty} \frac{u}{\sqrt{1+x^2}} = \lim_{x\rightarrow -\infty} \frac{u}{\sqrt{1+x^2}}$ exists, then
\begin{align*}
\lim_{M,N\rightarrow +\infty} \Big( \int_{-N}^{M} \Big| \frac{\ed u}{\ed x} \Big|^2 \ed x - \frac{x u^2}{1+x^2} \Big|_{-N}^{M} \Big)
\geq
-\frac{4\pi^2}{\int_{-\infty}^{+\infty}  u^{-2} \ed x},
\end{align*}
with the equality being achieved at $u^{-2} =  \frac{2k\sqrt{1-\alpha^2}}{1+ \alpha \big[\frac{1-x_0^2}{1+x_0^2}\cdot (1-x^2) + \frac{2x_0}{1+x_0^2}\cdot 2x \big]}$, where $k>0$, $\alpha \in (-1,1)$, $x_0 \in \mathbb{R}$.

\item If limits $\lim_{x\rightarrow +\infty} \frac{u}{\sqrt{1+x^2}} $, $\lim_{x\rightarrow -\infty} \frac{u}{\sqrt{1+x^2}} $ both exist, then
\begin{align*}
\lim_{M,N\rightarrow +\infty} 4 \Big( \int_{-N}^{M} \Big| \frac{\ed u}{\ed x} \Big|^2 \ed x - \frac{x u^2}{1+x^2} \Big|_{-N}^{M} \Big) + \int_{-\infty}^{+\infty}  \frac{3 u^2}{(1+x^2)^2} \ed x
\geq
-\frac{\pi^2}{\int_{-\infty}^{+\infty}  u^{-2} \ed x},
\end{align*}
with the equality being achieved at $u^{-2} =  \frac{2}{1+x^2} \cdot \frac{ k\sqrt{1-\alpha^2}}{1+ \frac{\alpha x}{\sqrt{1+x^2}}}$, where $k>0$, $\alpha \in (-1,1)$.
\end{enumerate}
\end{corollary}
\begin{proof}
From $u = \sqrt{\frac{1+x^2}{2}} v$, we obtain that
\begin{align*}
&
\frac{\ed u}{\ed x} =  \sqrt{\frac{1+x^2}{2}} \frac{\ed v}{\ed x} + \frac{1}{\sqrt{2}} \frac{x}{\sqrt{1+x^2}} v,
\\
&
\Big| \frac{\ed u}{\ed x} \Big|^2 = \frac{1+x^2}{2} \Big| \frac{\ed v}{\ed x} \Big|^2 + x v \cdot \frac{\ed v}{\ed x} + \frac{x^2}{2(1+x^2)} v^2,
\end{align*}
therefore
\begin{align*}
\int_{-\infty}^{+\infty} \Big| \frac{\ed u}{\ed x} \Big|^2 \ed x
&=
\int_{-\infty}^{+\infty}  \frac{1+x^2}{2} \Big| \frac{\ed v}{\ed x} \Big|^2 \ed x
+
\lim_{N\rightarrow +\infty} \int_{-N}^{N}  \Big[ x v \cdot \frac{\ed v}{\ed x} + \frac{x^2}{2(1+x^2)} v^2 \Big] \ed x
\\
&=
\int_{-\infty}^{+\infty}  \frac{1+x^2}{2} \Big| \frac{\ed v}{\ed x} \Big|^2 \ed x
-
\int_{-\infty}^{+\infty}  \frac{1}{2(1+x^2)} v^2 \ed x
+
\lim_{N\rightarrow +\infty} \frac{x v^2}{2}\Big|_{-N}^{N}
\\
&=
\int_{-\infty}^{+\infty}  \frac{1+x^2}{2} \Big| \frac{\ed v}{\ed x} \Big|^2 \ed x
-
\int_{-\infty}^{+\infty}  \frac{1}{2(1+x^2)} v^2 \ed x
+
\lim_{N\rightarrow +\infty} \frac{x u^2}{1+x^2} \Big|_{-N}^{N},
\end{align*}
thus $| \frac{\ed u}{\ed x} |^2$ is not integrable generally, and
\begin{align*}
\lim_{M,N\rightarrow +\infty} \Big( \int_{-N}^{M} \Big| \frac{\ed u}{\ed x} \Big|^2 \ed x - \frac{x u^2}{1+x^2} \Big|_{-N}^{M} \Big)
=
\int_{-\infty}^{+\infty}  \frac{1+x^2}{2} \Big| \frac{\ed v}{\ed x} \Big|^2 \ed x
-
\int_{-\infty}^{+\infty}  \frac{1}{2(1+x^2)} v^2 \ed x.
\end{align*}
Therefore the inequality can be rewritten in terms of $u$ as
\begin{align*}
\lim_{M,N\rightarrow +\infty} \Big( \int_{-N}^{M} \Big| \frac{\ed u}{\ed x} \Big|^2 \ed x - \frac{x u^2}{1+x^2} \Big|_{-N}^{M} \Big)
\geq
-\frac{\pi^2}{\int_{-\infty}^{+\infty}  u^{-2} \ed x},
\end{align*}
where $\lim_{x\rightarrow +\infty} \frac{u}{\sqrt{1+x^2}}= \lim_{x\rightarrow -\infty} \frac{u}{\sqrt{1+x^2}}$ exists. A similar derivation implies \textit{b.}
\end{proof}
\begin{remark}\label{rem 8.6}
Easy to see that similar inequalities as in corollary \ref{cor 8.4} hold.
\end{remark}

\section{Sketchy proof of the sharp Sobolev inequality on $\mathbb{S}^n, n\geq 3$}\label{sec 9}

In this section, we sketch another proof of the sharp Sobolev inequality \eqref{eqn 1.3} on $\mathbb{S}^n$, $n\geq 3$ with the method introduced in sections \ref{sec 6} and \ref{sec 7}. Recall the sharp Sobolev inequality \eqref{eqn 1.3} states that
\begin{align*}
\int_{\mathbb{S}^n} | \nabla v|^2 \ed \mu + \frac{n(n-2)}{4} \int_{\mathbb{S}^n} v^2 \ed \mu 
\geq
\frac{n(n-2)}{4} |\mathbb{S}^n|^{\frac{2}{n}} \Big( \int_{\mathbb{S}^n} v^{\frac{2n}{n-2}} \ed \mu \Big)^{\frac{n-2}{n}},
\quad
n\geq 3,
\end{align*}
and the equality is achieved at $v^{\frac{2n}{n-2}}(p)=\frac{k\sqrt{1-\alpha^2}}{1+ \alpha \cos d(p,p_0)}$ where $k>0$, $\alpha \in (-1,1)$, $p_0 \in \mathbb{S}^n$ and $d(p, p_0)$ is the distance between $p, p_0$.

We first introduce the relation between $\mathbb{S}^n$ and the Minkowski spacetime $\mathbb{M}^{n+2}$. Consider the coordinate system $\{ \theta, \vartheta\}$ of $\mathbb{S}^n$ where $\theta\in[0,\pi]$, $\vartheta \in \mathbb{S}^{n-1}$. Then we can introduce the spatial polar coordinate system $\{t, r, \theta, \vartheta\}$ of $\mathbb{M}^{n+2}$ where its transformation to the rectangular coordinate system is given by\footnote{We abuse the notation $\vartheta$ to denote the coordinate system of $\mathbb{S}^{n-1}$ and also the embedding of $\mathbb{S}^{n-1}$ into $\mathbb{E}^{n}$.}
\begin{align*}
(t, r, \theta, \vartheta) 
\mapsto
(t, r \cos \theta, r \sin \theta \cdot \vartheta).
\end{align*}
With the above $(t, r, \theta, \vartheta)$ coordinate system of $\mathbb{M}^{n+2}$, we introduce the Lorentz transformation $\varphi_{\alpha}$ similar to the one $\varphi_{\alpha, \theta_0, \btheta_0}$ in section \ref{sec 3} that
\begin{align*}
\varphi_{\alpha}:
\quad
(t, r, \theta, \vartheta) 
\quad
\mapsto 
\quad
(\bt, \br, \btheta, \bvartheta),
\end{align*}
where $\bt = \frac{t+\alpha r \cos \theta}{\sqrt{1-\alpha^2}}$, $\br=\sqrt{\frac{(r+\alpha t \cos \theta)^2 - \alpha^2(r^2 -t^2) \sin^2 \theta}{1-\alpha^2}}$, $\bvartheta=\vartheta$ and $\btheta$ is solved by
\begin{align*}
\textstyle
\cos \btheta 
=
\frac{\alpha t+ r\cos \theta}{\sqrt{(r+ \alpha t \cos \theta )^2 - \alpha^2 (r^2-t^2) \sin^2 \theta}},
\quad
\sin \btheta
=
\frac{r \sin \theta \sqrt{1-\alpha^2} }{ \sqrt{(r+ \alpha t \cos \theta )^2 - \alpha^2 (r^2-t^2) \sin^2 \theta} }.
\end{align*}
Restrict the Lorentz transformation $\varphi_{\alpha}$ on the lightcone $C_0$ with the coordinate system $(\uu, \theta, \vartheta)$ where $\uu=\frac{t+r}{2}$,
\begin{align*}
\textstyle
\varphi_{\alpha}:
\quad
(\uu, \theta ,\vartheta)
\quad
\mapsto
\quad
(\buu, \btheta, \bvartheta)
= 
( \frac{1+\alpha \cos \theta}{\sqrt{1-\alpha^2}} \uu, \btheta, \vartheta),
\end{align*}
where
\begin{align}
\textstyle
\cos \btheta 
=
\frac{\alpha + \cos \theta}{1+ \alpha  \cos \theta},
\quad
\sin \btheta
=
\frac{  \sqrt{1-\alpha^2} \sin \theta }{ 1+ \alpha  \cos \theta }.
\label{eqn 9.1}
\end{align}
Consider a section $S$ of $C_0$ which is parameterised as the graph of $\uu$ of a function $f$ over the $(\theta, \vartheta)$ domain in the $\{\uu, \theta, \vartheta\}$ coordinate system, then $\varphi_{\alpha}$ transforms $S$ to $\varphi_{\alpha}(S)$ parameterised by $\barf = \varphi_{\alpha}(f)$ that
\begin{align*}
\barf(\btheta, \bvartheta)
=
\frac{\sqrt{1-\alpha^2}}{1-\alpha\cos \btheta} f(\theta,\vartheta)
=
\frac{1+ \alpha \cos \theta}{\sqrt{1-\alpha^2}} f(\theta, \vartheta),
\end{align*}
where $\vartheta = \bvartheta$ and $\theta$ is determined by equation \eqref{eqn 9.1} of the Lorentz transformation $\varphi_{\alpha}$. Introduce $\nu_{\alpha} (\btheta) = \frac{\sqrt{1-\alpha^2}}{1-\alpha\cos \btheta}$, then the transformation $\varphi_{\alpha}$ of $f$ can be written as $\barf(\btheta, \bvartheta) = \nu_{\alpha}(\btheta) f(\theta, \vartheta)$. Let $f=v^{\frac{2}{n-2}}$, then we introduce the transformation $\gamma_{\alpha}$ for $v$ similar as the transformation $\gamma_{\alpha, \theta_0, \btheta_0}$ in proposition \ref{pro 3.1} that
\begin{align*}
\bv= \gamma_{\alpha} (v),
\quad
\bv(\btheta, \bvartheta) = [\nu_{\alpha}(\btheta)]^{\frac{n-2}{2}} v(\theta, \vartheta).
\end{align*}
Note that formally let $n=1$, the formulae above reduce to their corresponding ones in section \ref{sec 3}.

Introduce the functional $\calF[v]$ that
\begin{align*}
\calF[v] 
=
\int_{\mathbb{S}^n} | \nabla v|^2 \ed \mu
+
\frac{n(n-2)}{4} \int_{\mathbb{S}^n} v^2 \ed \mu,
\end{align*}
and the constraint $\int_{\mathbb{S}^n} v^{\frac{2n}{n-2}} \ed \mu = |\mathbb{S}^n| = \frac{2\pi^{\frac{n+1}{2}}}{\Gamma(\frac{n+1}{2})}$. Both the functional $F[v]$ and the constraint are invariant under the Lorentz transformation $\gamma_{\alpha}$.

The symmetric increasing rearrangement on the sphere is defined as follows: let $v$ be a function on $\mathbb{S}^n$ and $v^*$ be the symmetric increasing rearrangement of $v$ which is given by
\begin{align*}
v^*(\theta,\vartheta) = v^*(\theta) > 0,
\quad
v^*(\theta_1) \geq v^*(\theta_2) \text{ if } \theta_1 \geq \theta_2,
\quad
|L_v(t)| = |L_{v^*}(t)|,
\end{align*}
where $L_v(t) = \{ p \in \mathbb{S}^n: v(p) >t \}$ and similarly for $L_{v^*}(t)$. We have that
\begin{align*}
\Vert \nabla v \Vert_{L^2}  \geq \Vert \nabla v^* \Vert_{L^2},
\end{align*}
by the P\'olya-Szeg\H{o} inequality on the sphere. We introduce the following set of functions $S_{c}(\mathbb{S}^n)$ which is similar to $S_c$ in proposition \ref{pro 6.13}:
\begin{enumerate}[label=\textit{\alph*}.]
\item $v=v^* > 0$,
\item $v \in \mathrm{H}^1(\mathbb{S}^n)$,
\item $v$ is constant on $[\frac{\pi}{2}, \pi]$,
\item $\int_{\mathbb{S}^n} v^{\frac{2n}{n-2}} \ed \mu = |\mathbb{S}|^n$.
\end{enumerate}
Then we have the following analogy of proposition \ref{pro 6.13} for $\calF[v]$ on $\mathbb{S}^n$.
\begin{proposition}\label{pro 9.1}
We have that
\begin{align*}
\inf \Big\{ \calF[v] : v\in\mathrm{H}^1(\mathbb{S}^n), v>0, \int_{\mathbb{S}^n} v^{\frac{2n}{n-2}} \ed \mu = |\mathbb{S}|^n = 2\pi \Big\} = \inf_{v\in S_c(\mathbb{S}^n)} \{ \calF[v] \},
\end{align*}
and the infimum of $\calF$ on the left can be achieved if and only if that it can be achieved in $S_c(\mathbb{S}^n)$.
\end{proposition}
\begin{proof}[Sketch of proof]
The proof is parallel to the proof of proposition \ref{pro 6.13}. For any positive function $v\in \mathrm{H}^1(\mathbb{S}^n)$, we can construct a decreasing sequence $\{v_n\}_{n\in \mathbb{N}}$ of $\calF$ as in definition \ref{def 6.4} with $v^{\frac{2n}{n-2}}$ taking the role of $v^{-2}$. Then $\{\sup\{v_n\}\}$ is nonincreasing, $\int_{\mathbb{S}^n} v_n^{\frac{2n}{n-2}} \ed \mu$ remains constant and $\{\calF[v_n]\}$ is also nonincreasing, which is analogous to lemma \ref{lem 6.5}. The same argument as in the proofs of lemmas \ref{lem 6.10}, \ref{lem 6.11} implies that $\{v_n\}$ a uniform positive lower bound. Then there exists a subsequence limit of $\{v_n\}$, denoted by $v_{\infty}$ such that $v_{\infty} \in S_c(\mathbb{S}^n)$ and $\calF[v_{\infty}] \leq \calF[v]$. The proposition follows.
\end{proof}
\begin{remark}
In fact, if we release the condition \textit{a.} in $S_c(\mathbb{S}^n)$ by $v=v^* \geq 0$, then we actually donot need the uniform positive lower bound of $\{v_n\}$ in the proof. Then the proof of proposition \ref{pro 9.1} is much easier than proposition \ref{pro 6.13}, since the definition of $\{v_n\}$ on $\mathbb{S}^n, n\geq 3$ automatically gives the upper bound of $v_{\infty}$ while we need lemma \ref{lem 6.11} for $\{v_n\}$ on $\mathbb{S}^1$ for the upper bound. Another point which makes the case of $\mathbb{S}^n, n\geq 3$ easier is the positive sign of the coefficient of $\int_{\mathbb{S}^n} v^2\ed \mu$ in $\calF[v]$, opposite to the negative coefficient of $\int_{\mathbb{S}^1} v^2 \ed \theta$ in $\calF[v]$ on $\mathbb{S}^1$.
\end{remark}
Then by proposition \ref{pro 9.1}, in order to study the variational problem of the infimum of $\calF[v]$ is sufficient to study the infimum of $\calF[v]$ in $S_c(\mathbb{S}^n)$. It is easy to show that $\inf_{v\in S_c(\mathbb{S}^n)} \{ \calF[v] \}$ can be achieved in $S_c(\mathbb{S}^n)$ since $v\in S_c(\mathbb{S}^n)$ automatically has a upper bound $2^{\frac{n-2}{2n}}$ since
\begin{align*}
v^{\frac{2n}{n-2}} \leq v^{\frac{2n}{n-2}}(\pi) = \frac{\int_{\theta \in [\frac{\pi}{2},\pi]} v^{\frac{2n}{n-2}} \ed \mu}{|\mathbb{S}^n|/2} \leq \frac{\int_{\mathbb{S}^n} v^{\frac{2n}{n-2}} \ed \mu}{|\mathbb{S}^n|/2} =2.
\end{align*}
Then we see that any minimising sequence $\{v_n\}$ of $\calF[v]$ in $S_c(\mathbb{S}^n)$ has a subsequence limit $v_{\infty}$ such that $\calF[v_{\infty}]$ achieves $\inf_{v\in S_c(\mathbb{S}^n)} \{ \calF[v] \}$. It is easy to verify $v\in S_c(\mathbb{S}^n)$ except the assumption \textit{a.} that $v>0$. However we can show that $v\equiv 1$ by solving the Euler-Lagrange equation of $\calF$ explicitly since $v$ is just a function of the coordinate $\theta$. Then we prove the sharp Sobolev inequality \eqref{eqn 1.3} on $\mathbb{S}^n$, $n\geq 3$. We refer to \cite{Ta1976} for the equality case.

\section*{Acknowledgements}
\addcontentsline{toc}{section}{Acknowledgement}
The author acknowledges the support of the National Natural Science Foundation of China under Grant No. 12201338.

\appendix
\section{Proof of corollary \ref{cor 7.2}}\label{appen A}

It is sufficient to show the following claim.
\begin{claim}\label{cla A.1}
If $v\in \mathrm{H}^1[0,\pi]$ and $v(\pi)=0$, then $\int_0^{\pi} v^{-2} \ed \theta = +\infty$.
\end{claim}
In order to prove this claim, following the strategy stated in footnote \ref{footnote 4}, we study the minimum of the Dirichlet energy $E(v) = \int_0^{\pi} (v')^2 \ed \theta$ in the class of functions $B_{m,c}^M$ defined by
\begin{enumerate}[label=\alph*.]
\item $v\in \mathrm{H}^1([0,\pi])$,
\item $\min_{\theta\in [0,\pi]}\{ v(\theta) \} = m>0$, $\max_{\theta\in [0,\pi]}\{ v(\theta) \} =M>0$ and $\int_0^{\pi} v^{-2} \ed \theta = c \in (\frac{\pi}{M^2}, \frac{\pi}{m^2})$.
\end{enumerate}
We claim the following assertion which implies claim \ref{cla A.1}.
\begin{claim}\label{cla A.2}
The minimum of the Dirichlet energy $\min_{v\in B_{m,c}^M}\{E(v)\}$ as $m\rightarrow 0^+$ satisfies that
\begin{align*}
\min_{v\in B_{m,c}^M } \{ E(v) \} \rightarrow + \infty,
\text{ as }
m \rightarrow 0^+,
\end{align*}
where $c$, $M$, $M^{-1}$ could vary but are uniformly bounded when taking the limit.
\end{claim}
The claim \ref{cla A.2} follows from the explicit calculation of $\min_{v\in B_{m,c}^M}\{E(v)\}$. By the rearrangement inequality, we can restrict the class of functions $B_{m,c}^M$ to $B_{m,c}^{\prime M}$ by assuming the following additional condition,
\begin{enumerate}[label=\alph*.]\addtocounter{enumi}{2}
\item $v$ is monotonically nonincreasing, $v(\pi) = m$, $v(0) = M$.
\end{enumerate}
Then we have that $\min_{v\in B_{m,c}^M}\{E(v)\} = \min_{v\in B_{m,c}^{\prime M}}\{E(v)\}$.
The Euler-Lagrange equation of the Dirichlet energy $E(v)$ in $B_{m,c}^{\prime M}$ is
\begin{align*}
&
v'' = \lambda v^{-3}
\quad
\Leftrightarrow
\quad
v'^2 =  \mu - \lambda v^{-2}
\quad
\underset{v' \leq 0}{\Leftrightarrow}
\quad
\int^v \frac{\ed v}{\sqrt{\mu - \lambda v^{-2}}} = -\theta - \Lambda,
\\
\Rightarrow
&
\left\{
\begin{aligned}
&
\mu \neq 0:
\quad
\sqrt{\mu v^2 - \lambda} =-\mu(\theta + \Lambda)
\quad
\Rightarrow
\quad
v^2 = \mu (\theta+ \Lambda)^2 + \frac{\lambda}{\mu},
\\
&
\mu =0:
\quad
v^2 = -2\sqrt{-\lambda} (\theta+ \Lambda).
\end{aligned}
\right.
\end{align*}
We have the following lemma of the solution of the Euler-Lagrange equation of $E(v)$ in $B_{m,c}^{\prime M}$.

\begin{lemma}\label{lem A.3}
Let $M>m>0$ and $c\in (\frac{\pi}{M^2}, \frac{\pi}{m^2})$. The solution of the Euler-Lagrange equation of the Dirichlet energy $E=\int_0^{\pi} |v'|^2 \ed \theta$ in $B_{m,c}^{\prime M}$ is classified into several cases depending on the value of $c$.
\begin{enumerate}[label=\alph*.]
\item $c\in \Big(\frac{\pi}{M^2}, \frac{\pi}{2M\sqrt{M^2-m^2}} \log \big( \frac{M + \sqrt{M^2-m^2}}{M - \sqrt{M^2-m^2}} \big) \Big)$: the solution takes the following form that $\alpha$ is a parameter in $(0,\pi)$ and
\begin{align*}
v_0^2(\theta) = 
\left\{
\begin{aligned}
&
M^2,
\quad
\theta \in [0,\alpha],
\\
&
-\frac{M^2-m^2}{(\pi-\alpha)^2} (\theta-\alpha)^2 + M^2,
\quad
\theta \in (\alpha,\pi].
\end{aligned}
\right.
\end{align*}

\item $c \in \Big[\frac{\pi}{2M\sqrt{M^2-m^2}} \log \big( \frac{M + \sqrt{M^2-m^2}}{M - \sqrt{M^2-m^2}} \big), \frac{\pi}{M^2-m^2} \log \frac{M^2}{m^2} \Big)$: the solution takes the following form 
\begin{align*}
v_0^2(\theta) = \mu (\theta+ \Lambda)^2 + \frac{\lambda}{\mu}.
\end{align*}
where $\mu<0$, $\lambda$, $\Lambda$ are parameters which can be determined by $\frac{\lambda}{\mu} \in [M^2, +\infty)$, $m$, $M$ from equations $v_0(\pi) = m$, $v_0(0) = M$.

\item $c=\frac{\pi}{M^2-m^2} \log \frac{M^2}{m^2}$: the solution is given by
\begin{align*}
v_0^2(\theta) = -\frac{M^2 - m^2}{\pi} \theta +M^2.
\end{align*}

\item $c \in \Big(\frac{\pi}{M^2-m^2} \log \frac{M^2}{m^2}, \frac{\pi}{m\sqrt{M^2-m^2}} \arctan \frac{\sqrt{M^2-m^2}}{m} \Big]$: the solution takes the following form
\begin{align*}
v_0^2(\theta) = \mu (\theta+ \Lambda)^2 + \frac{\lambda}{\mu},
\end{align*}
where $\mu>0$, $\lambda$, $\Lambda$ are parameters which can be determined by $\frac{\lambda}{\mu} \in (-\infty, m^2]$, $m$, $M$ from equations $v_0(\pi) = m$, $v_0(0) = M$.

\item $c \in \Big( \frac{\pi}{m\sqrt{M^2-m^2}} \arctan \frac{\sqrt{M^2-m^2}}{m}, \frac{\pi}{m^2} \Big)$: the solution takes the following form that $\beta$ is a parameter in $(0,\pi)$ and
\begin{align*}
v_0^2(\theta) = 
\left\{
\begin{aligned}
&
\frac{M^2-m^2}{\beta^2} (\theta-\beta)^2 + m^2,
\quad
\theta \in [0,\beta),
\\
&
m^2,
\quad
\theta \in [\beta,\pi].
\end{aligned}
\right.
\end{align*}
\end{enumerate}
The above solution $v_0$ achieves the minimum of the Dirichlet energy $\min_{v\in B_{m,c}^{\prime M}}\{E(v)\}$. The solutions in above cases are illustrated in figures \ref{fig 10} and \ref{fig 11}.

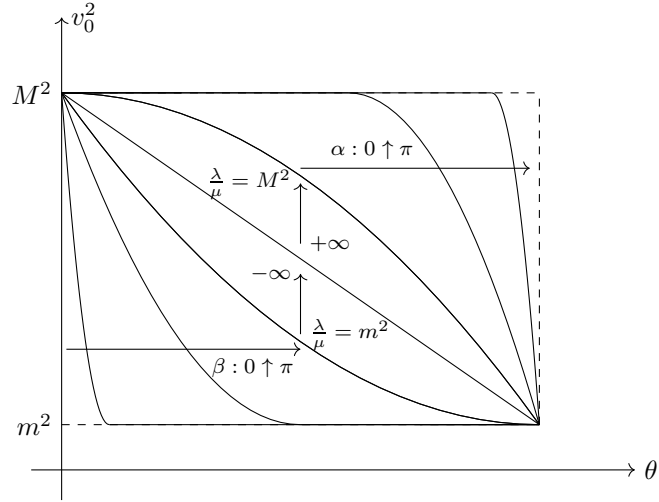
\begin{figure}[H]
\centering
\begin{tikzpicture}[scale=2]
\draw[->] (0,-0.2) -- (0,3) node[right] {$v_0^2$};
\draw[->] (-0.2,0) -- (1.2*pi,0) node[right] {$\theta$};
\node[left] at (0,0.3) {$m^2$};
\node[left] at (0,2.5) {$M^2$};
\draw[domain=0:pi] plot (\x, {-(2.5-0.3)*(\x)^2/(pi*pi)+2.5});
\draw[domain=0:pi/4] plot (\x, {2.5});
\draw[<-] (pi/2, 1.9) node[left] {\footnotesize $\frac{\lambda}{\mu}=M^2$} -- (pi/2,1.5) node[right] {\footnotesize $+\infty$};
\draw[domain=0:pi] plot (\x, {(2.5-0.3)*(\x-pi)^2/(pi*pi)+0.3});
\draw[domain=3*pi/4:pi] plot (\x, {0.3});
\draw[->] (pi/2,.9) node[right] {\footnotesize $\frac{\lambda}{\mu}=m^2$} -- (pi/2,1.3) node[left] {\footnotesize $-\infty$};
\draw[domain=0:pi] plot (\x, {-(2.5-0.3)*\x/pi+2.5});
\draw[domain=0:pi] plot (\x, {(2.5-0.3)*(\x-pi)^2/(pi*pi)+0.3});
\draw[domain=0:0.5*pi] plot (\x, {(2.5-0.3)*(\x-0.5*pi)^2/(0.5*pi*0.5*pi)+0.3});
\draw[domain=0.5*pi:pi] plot (\x, {0.3});
\draw[domain=0:0.1*pi] plot (\x, {(2.5-0.3)*(\x-0.1*pi)^2/(0.1*pi*0.1*pi)+0.3});
\draw[domain=0.1*pi:pi] plot (\x, {0.3});
\draw[->] (0.01*pi,0.8) -- (0.4*pi,0.8) node[below] {\footnotesize $\beta: 0 \uparrow \pi$} --  (0.5*pi,0.8);
\draw[dashed,domain=0:pi] plot (\x, {0.3});
\draw[domain=0:pi] plot (\x, {-(2.5-0.3)*(\x)^2/(pi*pi)+2.5});
\draw[domain=0:0.6*pi] plot (\x, {2.5});
\draw[domain=0.6*pi:pi] plot (\x, {-(2.5-0.3)*(\x-0.6*pi)^2/(0.4*pi*0.4*pi)+2.5});
\draw[domain=0:0.9*pi] plot (\x, {2.5});
\draw[domain=0.9*pi:pi] plot (\x, {-(2.5-0.3)*(\x-0.9*pi)^2/(0.1*pi*0.1*pi)+2.5});
\draw[dashed,domain=0:pi] plot (\x, {2.5});
\draw[dashed] (pi,0.3) -- (pi,2.5);
\draw[->] (0.5*pi,2) -- (0.65*pi,2) node[above] {\footnotesize $\alpha: 0 \uparrow \pi$} --  (0.98*pi,2);
\end{tikzpicture}
\caption{Summary of $v_0$. The parameter $c=\int_0^{\pi} v_0^{-2} \ed \theta$ decreases along the direction of the arrow.} 
\label{fig 10}
\end{figure}

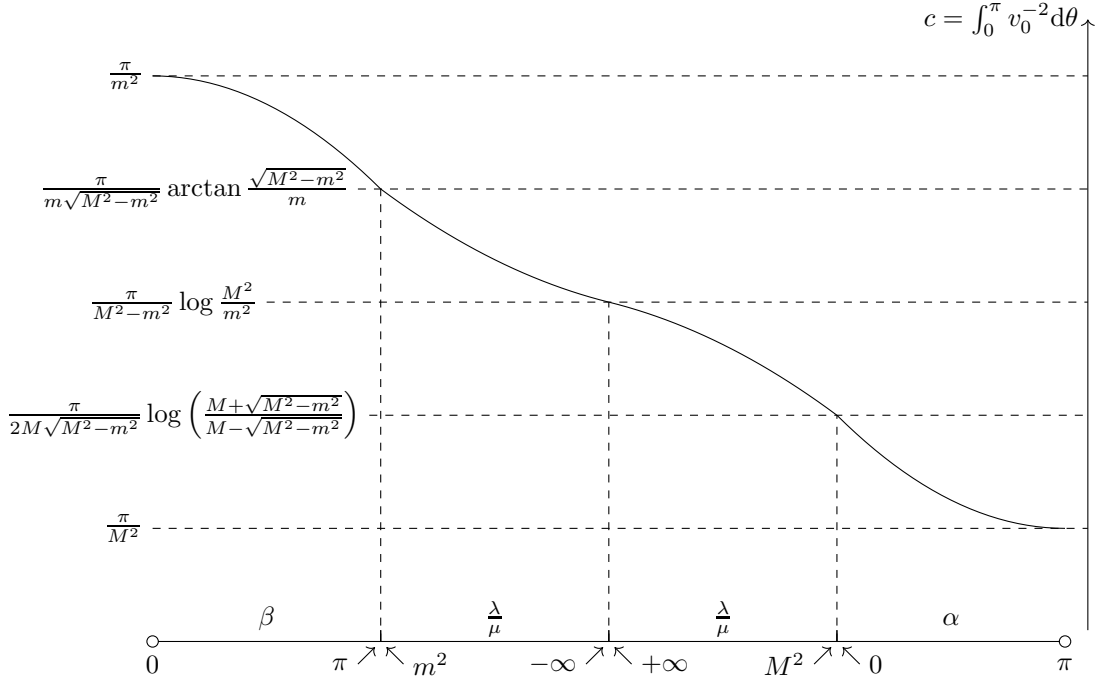
\begin{figure}[H]
\centering
\begin{tikzpicture}[scale=1.5]
\draw (0,0) circle(.05) (.05,0) -- (2,0) -- (2, .1);
\node at (0,-.2) {$0$};
\draw[->] (2-.2,-.2) node[left]{$\pi$} -- (2-.05,-.05);
\node at (1,.2) {$\beta$};
\draw (2, .1) -- (2,0) -- (4,0) -- (4, .1);
\draw[->] (2+.2,-.2) node[right]{$m^2$} -- (2+.05,-.05);
\draw[->] (4-.2,-.2) node[left]{$-\infty$} -- (4-.05,-.05);
\node at (3,.2) {$\frac{\lambda}{\mu}$};
\draw (4, .1) -- (4,0) -- (6,0) -- (6, .1);
\draw[->] (4+.2,-.2) node[right]{$+\infty$} -- (4+.05,-.05);
\draw[->] (6-.2,-.2) node[left]{$M^2$} -- (6-.05,-.05);
\node at (5,.2) {$\frac{\lambda}{\mu}$};
\draw (6, .1) -- (6,0) -- (8-.05,0)  (8,0) circle(.05);
\node at (8,-.2) {$\pi$};
\draw[->] (6+.2,-.2) node[right]{$0$} -- (6+.05,-.05);
\node at (7,.2) {$\alpha$};
\draw[dashed] (0,1) node[left] {$\frac{\pi}{M^2}$} -- (8.2,1) 
(2-.1,2) node[left] {$\frac{\pi}{2M\sqrt{M^2-m^2}} \log \Big( \frac{M + \sqrt{M^2-m^2}}{M - \sqrt{M^2-m^2}} \Big)$} -- (8.2,2) 
(1,3) node[left] {$\frac{\pi}{M^2-m^2} \log \frac{M^2}{m^2}$} -- (8.2,3) 
(1.8,4) node[left] {$\frac{\pi}{m\sqrt{M^2-m^2}} \arctan \frac{\sqrt{M^2-m^2}}{m}$} -- (8.2,4)
(0,5) node[left] {$\frac{\pi}{m^2}$} -- (8.2,5);
\draw[domain=0:2] plot (\x, {-\x*\x/(2*2)+5});
\draw[domain=2:4] plot (\x, {(\x-5)*(\x-5)/8 - 1/8 + 3});
\draw[domain=4:6] plot (\x, {-(\x-3)*(\x-3)/8 + 1/8 + 3});
\draw[domain=6:8] plot (\x, {(\x-8)*(\x-8)/(2*2)+1});
\draw[dashed] (2,0) -- (2,4);
\draw[dashed] (4,0) -- (4,3);
\draw[dashed] (6,0) -- (6,2);
\draw[->] (8.2,.1) -- (8.2,5.5) node[left] {$c=\int_0^{\pi} v_0^{-2} \ed \theta$};
\end{tikzpicture}
\caption{Summary of the parameter $c=\int_0^{\pi} v_0^{-2} \ed \theta$.}
\label{fig 11}
\end{figure}
\end{lemma}

Then we study the dependence of the minimum of the Dirichlet energy $\min_{v\in B_{m,c}^{\prime M}} \{E(v)\}$ on the parameter $c= \int_0^{\pi} v^{-2} \ed \theta$. We prove a lemma on the monotonicity of $\min_{v\in B_{m,c}^{\prime M}} \{E(v)\}$ with respect to the parameter $c$ first.

\begin{lemma}\label{lem A.4}
The derivative $\frac{\ed}{\ed c} E(v_0)$ of the minimum of the Dirichlet energy $E$ in $B_{m,c}^{\prime M}$ with respect to $c$ is $\lambda$, i.e. $\frac{\ed}{\ed c} E(v_0)=\lambda$. Thus the monotonicity of $E(v_0)$ with respect to $c$ is determined by the sign fo $\lambda$.
\end{lemma}
\begin{proof}
Calculate the derivative $\frac{\ed}{\ed c} E(v_0)$. Use the dot on the top $\cdot$ to denote the derivative $\frac{\ed}{\ed c}$. Let $\beta= \pi$ in cases \textit{a. b. c. d.}, and $\alpha=0$ in cases \textit{b. c. d. e.} in lemma \ref{lem A.3}. Then
\begin{align*}
\frac{\ed}{\ed c} E(v_0)
&=
\frac{\ed}{\ed c} \int_{\alpha}^{\beta} |v_0'|^2 \ed \theta
\\
&= 
\dot{\beta} |v_0'(\beta) |^2 
- \dot{\alpha} |v_0'(\alpha)|^2 
+ 2 \int_{\alpha}^{\beta} v_0' \cdot \dot{v}_0' \ed \theta
=
- 2 \int_{\alpha}^{\beta} v_0'' \cdot \dot{v}_0 \ed \theta 
+ 2v_0' \cdot \dot{v}_0 \big\vert_{\alpha}^{\beta}
\\
&=
- 2 \int_{\alpha}^{\beta} \lambda v_0^{-3} \cdot \dot{v}_0 \ed \theta 
=
\lambda \int_{\alpha}^{\beta} \frac{\ed}{\ed c}(v_0^{-2}) \ed \theta 
=
\lambda \int_0^{\pi} \frac{\ed}{\ed c}(v_0^{-2}) \ed \theta
\\&=
\lambda \frac{\ed}{\ed c} \int_0^{\pi} v_0^{-2} \ed \theta
=
\lambda.
\end{align*}
\end{proof}

Then we can obtain the detailed dependence of $\min_{v\in B_{m,c}^{\prime M}} \{E(v)\}$ on the parameter $c= \int_0^{\pi} v^{-2} \ed \theta$.
\begin{lemma}\label{lem A.5}
Adopt the notations and the classification of the minimiser $v_0$ of the Dirichlet energy $E(v)$ in $B_{m,c}^{\prime M}$ in lemma \ref{lem A.3}.
\begin{enumerate}
\item
There is a unique value of $c$, denoted by $c_{\lambda=0}$, solving $\frac{\ed}{\ed c} E(v_0) = \lambda = 0$. $c_{\lambda=0}$ lies in the interval $\Big(\frac{\pi}{M^2-m^2} \log \frac{M^2}{m^2}, \frac{\pi}{m\sqrt{M^2-m^2}} \arctan \frac{\sqrt{M^2-m^2}}{m} \Big]$, which corresponds to case d. in lemma \ref{lem A.3}.

\item
$\frac{\ed}{\ed c} E(v_0) = \lambda<0$ for $c \in (\frac{\pi}{M^2}, c_{\lambda=0})$ and $\frac{\ed}{\ed c} E(v_0) = \lambda>0$ for $c\in (c_{\lambda=0}, \frac{\pi}{m^2})$.

\item
$\lambda \rightarrow - \infty$ as $c \rightarrow (\frac{\pi}{M^2})^+$ and $\lambda \rightarrow + \infty$ as $c \rightarrow (\frac{\pi}{m^2})^-$ .

\item
The minimiser at $c_{\lambda=0}$ is denoted by $v_{0,c_{\lambda=0}}$. We have that
\begin{align*}
v_{0,c_{\lambda=0}}^2(\theta) 
=
\mu (\theta+\Lambda)^2
=
\frac{(M-m)^2}{\pi^2} \big( \theta - \frac{\pi M}{M-m} \big)^2.
\end{align*}

\item $c_{\lambda=0} = \frac{\pi}{mM}$.

\item The Dirichlet energy $E(v_{0,c_{\lambda=0}}) = \frac{(M-m)^2}{\pi}$.

\item Denote the minimiser of $E$ in $B_{m,c}^M$ by $v_{0,c}$, then 
\begin{align*}
\lim_{c\rightarrow (\frac{\pi}{M^2})^+} E(v_{0,c}) = \lim_{c\rightarrow (\frac{\pi}{m^2})^-} E(v_{0,c}) = +\infty.
\end{align*}
\end{enumerate}
See the illustration of the lemma in figure \ref{fig 12}.
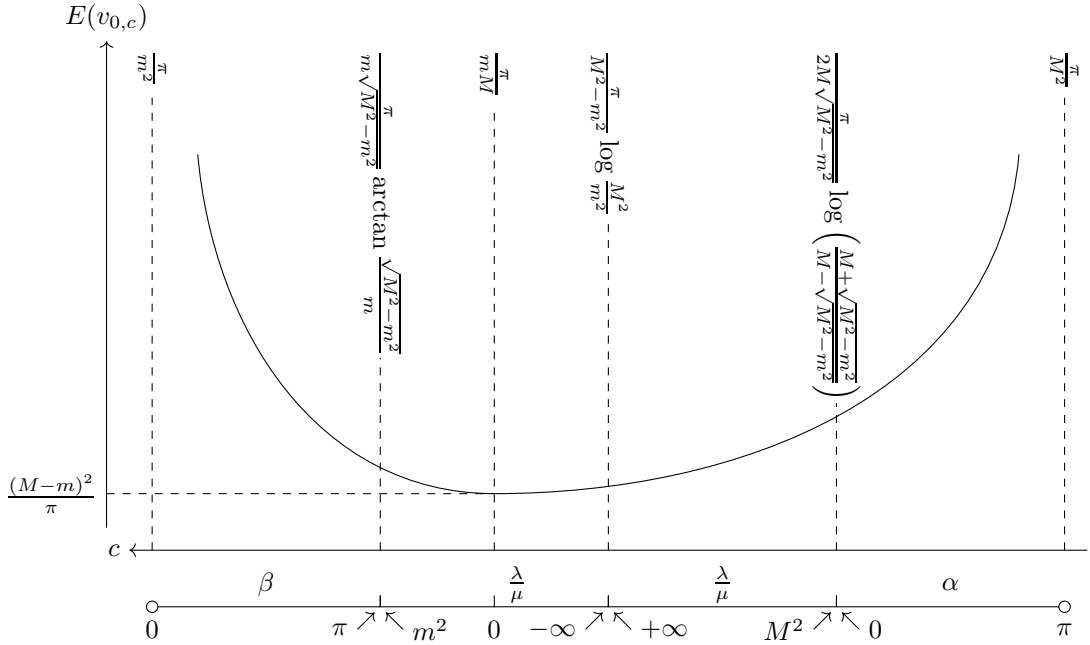
\begin{figure}[H]
\centering
\begin{tikzpicture}[scale=1.5]
\draw (0,0) circle(.05) (.05,0) -- (2,0) -- (2, .1);
\node at (0,-.2) {$0$};
\draw[->] (2-.2,-.2) node[left]{$\pi$} -- (2-.05,-.05);
\node at (1,.2) {$\beta$};
\draw (2, .1) -- (2,0) -- (3,0) -- (3,.1) (3,0) -- (4,0) -- (4, .1);
\node at (3,-.2) {$0$};
\draw[->] (2+.2,-.2) node[right]{$m^2$} -- (2+.05,-.05);
\draw[->] (4-.2,-.2) node[left]{$-\infty$} -- (4-.05,-.05);
\node at (3.2,.2) {$\frac{\lambda}{\mu}$};
\draw (4, .1) -- (4,0) -- (6,0) -- (6, .1);
\draw[->] (4+.2,-.2) node[right]{$+\infty$} -- (4+.05,-.05);
\draw[->] (6-.2,-.2) node[left]{$M^2$} -- (6-.05,-.05);
\node at (5,.2) {$\frac{\lambda}{\mu}$};
\draw (6, .1) -- (6,0) -- (8-.05,0)  (8,0) circle(.05);
\node at (8,-.2) {$\pi$};
\draw[->] (6+.2,-.2) node[right]{$0$} -- (6+.05,-.05);
\node at (7,.2) {$\alpha$};
\draw[<-] (-.2,.5) node[left] {$c$} -- (0,.5) -- (2,.5) -- (4,.5) -- (6,.5) -- (8,.5) -- (8.2,.5);
\draw[dashed] (0,.5) -- (0,5-.5)
(2,.5) -- (2,5-2.8)
(3,.5) -- (3,5-.6)
(4,.5) -- (4,5-1.6)
(6,.5) -- (6,5-3.2)
(8,.5) -- (8,5-.5);
\node[rotate=-90,right] at (0,5) {$\frac{\pi}{m^2}$};
\node[rotate=-90,right] at (2,5) {$\frac{\pi}{m\sqrt{M^2-m^2}} \arctan \frac{\sqrt{M^2-m^2}}{m}$};
\node[rotate=-90,right] at (3,5) {$\frac{\pi}{mM}$};
\node[rotate=-90,right] at (4,5) {$\frac{\pi}{M^2-m^2} \log \frac{M^2}{m^2}$};
\node[rotate=-90,right] at (6,5) {$\frac{\pi}{2M\sqrt{M^2-m^2}} \log \Big( \frac{M + \sqrt{M^2-m^2}}{M - \sqrt{M^2-m^2}} \Big)$};
\node[rotate=-90,right] at (8,5) {$\frac{\pi}{M^2}$};
\draw[->] (-.4,.7) -- (-.4,5) node[above] {$E(v_{0,c})$};
\draw (.4,4)
[out=-85,in=180] to (3,1)
[out=0,in=-95] to (8-.4,4);
\draw[dashed] (-.4,1) node[left] {$\frac{(M-m)^2}{\pi}$} -- (3,1);
\end{tikzpicture}
\caption{Summary of the minimum of the Dirichlet energy in $B_{m,c}^M$: $E(v_{0,c})$.}
\label{fig 12}
\end{figure}
\end{lemma}
\begin{proof}
All can be proved by explicit calculations. We just give a brief proof of \textit{7.}. In case \textit{a.} of lemma \ref{lem A.3}, we have that
\begin{align}
\begin{aligned}
E(v_{0,c}) 
&=
\mu(\pi -\alpha) - \lambda \int_{\alpha}^{\pi} v_{0,c}^{-2} \ed \theta
\\
&=
-\frac{M^2-m^2}{\pi-\alpha} + \frac{M\sqrt{M^2-m^2}}{2(\pi-\alpha)} \log \Big( \frac{M + \sqrt{M^2-m^2}}{M-\sqrt{M^2 -m^2}}  \Big)
\rightarrow +\infty,
\quad \alpha \rightarrow \pi^-.
\end{aligned}
\label{eqn A.1}
\end{align}
In case \textit{e.} of lemma \ref{lem A.3}, we have that
\begin{align*}
E(v_{0,c}) 
&=
\mu \beta - \lambda \int_0^{\beta} v_{0,c}^{-2} \ed \theta
\\
&=
\frac{M^2-m^2}{\beta} - \frac{m\sqrt{M^2-m^2}}{\beta} \arctan \frac{\sqrt{M^2-m^2}}{m}
\rightarrow +\infty,
\quad \beta \rightarrow 0^+.
\end{align*}
Then \textit{7.} follows.
\end{proof}

Now we are ready to prove claim \ref{cla A.2}. It follows from the following lemma.
\begin{lemma}\label{lem A.6}
Given $M$, $M^{-1}$ and $c > \pi/M^2$ all uniformly bounded during taking the limit, the minimum of the Dirichlet energy $E$ in $B_{m,c}^{\prime M}$ approaches $+\infty$ as $m \rightarrow 0^+$.
\end{lemma}
In order to prove the above lemma, we introduce some notations first. See figure \ref{fig 13}.
\begin{enumerate}[label=\textbullet]
\item $c_{\alpha=0}$: the corresponding value of $c$ in case \emph{a.} of lemma \ref{lem A.3} when $\alpha=0$. We have that
\begin{align*}
c_{\alpha=0}
=
\frac{\pi}{2M\sqrt{M^2-m^2}} \log \Big( \frac{M + \sqrt{M^2-m^2}}{M - \sqrt{M^2-m^2}} \Big).
\end{align*}

\item $v_{0, c_{\alpha=0}}$: the minimiser of $E$ in $B_{m,c_{\alpha=0}}^M$. We have that as obtained in lemma \ref{lem A.3}
\begin{align*}
v_{0,c_{\alpha=0}}^2(\theta) = -\frac{M^2-m^2}{\pi^2} \theta^2 +M^2.
\end{align*}

\item $E(v_{0, c_{\alpha=0}})$:  the Dirichlet energy of $v_{0,c_{\alpha=0}}$, which is the minimum of $E$ in $B_{m,c_{\alpha=0}}^M$. By equation \eqref{eqn A.1},
\begin{align*}
E(v_{0,c_{\alpha=0}})
=
-\frac{M^2-m^2}{\pi} + \frac{M\sqrt{M^2-m^2}}{2\pi} \log \Big( \frac{M + \sqrt{M^2-m^2}}{M-\sqrt{M^2 -m^2}}  \Big).
\end{align*}
\end{enumerate}
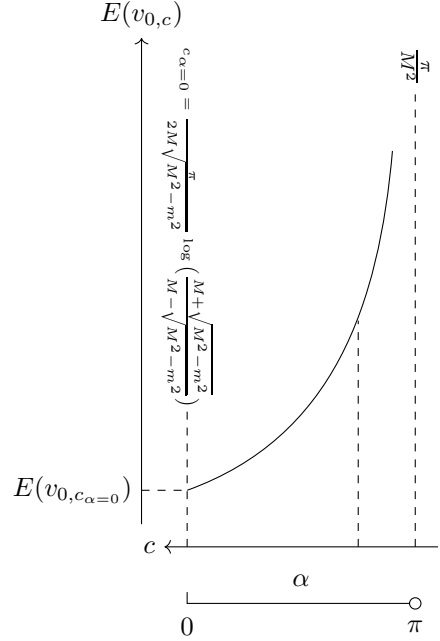
\begin{figure}[H]
\centering
\begin{tikzpicture}[scale=1.5]
\draw (6, .1) -- (6,0) -- (8-.05,0)  (8,0) circle(.05);
\node at (8,-.2) {$\pi$};
\node at (6,-.2) {$0$};
\node at (7,.2) {$\alpha$};
\draw[<-] (6-.2,.5) node[left] {$c$} -- (6,.5) -- (8,.5) -- (8.2,.5);
\draw[dashed] 
(6,.5) -- (6,5-3.3)
(8,.5) -- (8,5-.5);
\node[rotate=-90,right] at (6,5) {\tiny $c_{\alpha=0}=\frac{\pi}{2M\sqrt{M^2-m^2}} \log \Big( \frac{M + \sqrt{M^2-m^2}}{M - \sqrt{M^2-m^2}} \Big)$};
\node[rotate=-90,right] at (8,5) {$\frac{\pi}{M^2}$};
\draw[->] (6-.4,.7) -- (6-.4,5) node[above] {$E(v_{0,c})$};
\draw (6,1)
[out=20,in=-95] to (8-.2,4);
\draw[dashed] (6-.4,1) node[left] {$E(v_{0,c_{\alpha=0}})$} -- (6,1);
\draw[dashed] 
(7.5,.5) -- (7.5,5-2.5);
\end{tikzpicture}
\caption{$c_{\alpha=0}$ and $E(v_{0,c_{\alpha=0}})$.}
\label{fig 13}
\end{figure}
A simple calculation shows that $c_{\alpha=0} \rightarrow +\infty$, $E(v_{0,c_{\alpha=0}}) \rightarrow +\infty$ as $m \rightarrow 0^+$ given $M$, $M^{-1}$ both bounded, thus we have $E(v_{0,c}) \rightarrow +\infty$ as $m\rightarrow^+$ given $M$, $M^{-1}$ and $c > \frac{\pi}{M^2}$ all bounded.
\begin{proof}[Proof of lemma \ref{lem A.6}]
Since $c_{\alpha=0} \rightarrow +\infty$ as $m\rightarrow 0^+$, then $c\in (\frac{\pi}{M^2}, c_{\alpha=0})$ when $m$ is sufficiently small. Thus
\begin{align*}
E(v_{0,c}) > E (v_{0,c_{\alpha=0}}) \rightarrow + \infty,
\text{ as } m\rightarrow 0^+,
\end{align*}
since $\frac{\ed}{\ed c} E(v_{0,c}) <0$ for any $c\in (\frac{\pi}{M^2}, c_{\alpha=0})$ by lemma \ref{lem A.5}. See figure \ref{fig 13} for the illustration of the proof.
\end{proof}

\Address

\end{document}